\documentclass[12pt]{amsart}
\usepackage{latexsym}		%
\usepackage{amssymb}
\usepackage{epsfig}		%
\usepackage{amsfonts}

\renewcommand{\baselinestretch}{1.1}		%
\normalsize					%

\newcommand{\IA}{\mathbb{A}}

\newcommand{\IP}{\mathbb{P}}                                     
\newcommand{\IR}{\mathbb{R}}                           
\newcommand{\IC}{\mathbb{C}}
\newcommand{\IZ}{\mathbb{Z}}

\newcommand{\M}{\mathcal{M}}
\newcommand{\Mb}{\mathcal{M}'}
\newcommand{\K}{\mathcal{K}}

\newcommand{\R}{\mathcal{R}}

\newcommand{\cO}{\mathcal{O}}

\newcommand{\G}{\mathcal{G}}
\newcommand{\cH}{\mathcal{H}}

\newcommand{\cS}{\mathcal{S}}

\newcommand{\g}{       \mathfrak{g}     }

\newcommand{\gl}{       \mathfrak{gl}     } %
\newcommand{\so}{       \mathfrak{so}     }

\newcommand{\lk}{\mathfrak{k}}

\newcommand{\lsl}{\mathfrak{sl}} 
\newcommand{\Lia}{{\text{\rm Lie}}}	%

\newcommand{\lt}{\mathfrak{t}}

\newcommand{\ld}{\mathfrak{d}}

\newcommand{\disc}{\Delta}
\newcommand{\Ds}{\Delta^{\!*}}

\newcommand{\pf}{\begin{bpf}}

\newcommand{\pfms}{\begin{bpfms}}
\newcommand{\epf}{\end{bpf}\hfill$\square$\\}		%
\newcommand{\epfms}{\end{bpfms}\hfill$\square$\\}		%

\newcommand{\lu}{\mathfrak{u}}
\newcommand{\wt}{\widetilde}

\newcommand{\wh}{\widehat}
\newcommand{\al}{\alpha}
\newcommand{\be}{\beta}
\newcommand{\ga}{\gamma}
\newcommand{\Ga}{\Gamma}
\newcommand{\si}{\sigma}
\newcommand{\Si}{\Sigma}
\newcommand{\Sif}{\Sigma_i(\wh F)}

\newcommand{\bd}{{\bf d}}
\newcommand{\bdp}{{\bf d'}}

\newcommand{\bk}{{\bf k}}

\newcommand{\ISto}{\text{\rm $\mathbb{S}$to}}

\newcommand{\Sect}{\text{\rm Sect}}
\newcommand{\Ssect}{\wh {\text{\rm Sect}}}

\newcommand{\Ad}{\text{\rm Ad}}
\newcommand{\ad}{\text{\rm ad}}

\newcommand{\pr}{\text{\rm pr}}		%

\newcommand{\tr}{\text{\rm Tr}}
\newcommand{\Hom}{\text{\rm Hom}}

\newcommand{\End}{\text{\rm End}}
\newcommand{\diag}{{\text{\rm diag}}}

\newcommand{\od}{{\text{\footnotesize\rm od}}}		%

\newcommand{\reg}{{\text{\footnotesize\rm reg}}}
\newcommand{\Real}{{\text{\rm Re}}}

\newcommand {\flb}{\lbrack\!\lbrack}
\newcommand {\frb}{\rbrack\!\rbrack}

\newcommand {\pibpm}{{^i\!\,b_\pm}}
\newcommand {\pibm}{{^i\!\,b_-}}
\newcommand {\pibp}{{^i\!\,b_+}}

\def\mapright#1{\smash{
	\mathop{\longrightarrow}\limits^{#1}}}

\def\mapdown#1{\Big\downarrow
	\rlap{$\vcenter{\hbox{$\scriptstyle#1$}}$}}

\theoremstyle{plain}
\newtheorem {hypo}{\bf\hspace{-\parindent}Hypothesis}
\newtheorem*{plainthm}{Theorem}
\newtheorem {thm}{Theorem}
\newtheorem {prop}[hypo]{Proposition}%

\newtheorem {cor}[hypo]{Corollary}%
\newtheorem {lem}[hypo]{Lemma}%

\theoremstyle{definition}
\newtheorem {defn}[hypo]{Definition}%
\newtheorem {eg}[hypo]{Example}%

\theoremstyle{remark}
\newtheorem {rmk}[hypo]{Remark}%

\begin{document}

\title[$G$-bundles, isomonodromy and quantum Weyl groups]{$\bf G$-bundles, 
isomonodromy and \\
quantum Weyl groups}
\author{P. P. Boalch}
\address{I.R.M.A., Universit\'e Louis Pasteur et C.N.R.S.\\
7 rue Ren\'e Descartes\\
67084 Strasbourg Cedex\\
France}
\email{boalch@math.u-strasbg.fr}

\maketitle

\renewcommand{\baselinestretch}{1.1}		%
\normalsize

\begin{section}{Introduction}

It is now twenty years since Jimbo, Miwa and Ueno \cite{JMU81}
generalised Schlesinger's equations (governing isomonodromic deformations 
of logarithmic connections on vector bundles over the Riemann sphere) 
to the case of connections with arbitrary order poles.
An interesting feature was that new deformation parameters arose:
one may vary the `irregular type' of the connections at each
pole of order two or more (irregular pole), as well as the pole positions.
Indeed, for each irregular pole the fundamental group of the space of 
deformation parameters was multiplied by a factor of
$$P_n= \pi_1\left(\IC^n\setminus \text{diagonals}\right),$$
where $n$ is the rank of the vector bundles.
(This factor arose because the connections must be `generic'; the leading
term at each irregular pole must have distinct eigenvalues.)

The motivation behind the first part of this paper is
the question of
how to generalise the work of Jimbo, Miwa and Ueno 
(and also \cite{smpgafm, smid})
to the case of meromorphic connections on principal $G$-bundles
for complex reductive groups $G$. 
For simple poles (Schlesinger's equations) this generalisation is
immediate, but in general one needs to understand the `$G$-valued' Stokes
phenomenon in order to proceed 
(that is, one needs   to understand the local moduli of meromorphic
connections on $G$-bundles). 
This will be done in Section \ref{sn: sms}. 
Naturally enough a good theory is obtained provided the leading
term at each irregular pole is {\em regular semisimple} (that is, 
lies on the complement  of the root hyperplanes in some Cartan subalgebra).
The main result of Section \ref{sn: sms} is an irregular Riemann--Hilbert
correspondence describing the local moduli in terms
of $G$-valued Stokes multipliers, and is the natural generalisation of 
the result of Balser, Jurkat and Lutz \cite{BJL79} in the $GL_n(\IC)$ case.
The proof is necessarily quite different to that of \cite{BJL79}
however.

In the rest of the paper we consider isomonodromic deformations of such
connections in the simplest case: that of connections with one order two pole
over the unit disc.
The main things we will prove are: 
1) That the classical actions of
quantum Weyl groups found by De \!Concini, Kac and Procesi \cite{DKP}
do arise from isomonodromy (and so have a purely geometrical origin)
and 2) That a certain flat connection appearing in work of 
De \!Concini and Toledano Laredo arises directly from the isomonodromy
Hamiltonians,
indicating that the previous result is the classical analogue of their
conjectural Kohno--Drinfeld theorem for quantum Weyl groups.

In more detail, in this `simplest case' 
the fundamental group of the space of deformation parameters is the 
generalised pure braid group associated to $\g=\Lia(G)$:
$$P_\g=\pi_1(\lt_\reg)$$
where $\lt_\reg$ is the regular subset of a Cartan subalgebra $\lt\subset \g$.
By considering isomonodromic deformations one obtains a nonlinear 
(Poisson) action of $P_\g$ as follows
(this is purely geometrical---as explained in \cite{smid} 
the author likes to think of isomonodromy as
a natural analogue of the Gauss--Manin connection in non-Abelian cohomology):
There is a moduli space $\M$ of generic (compatibly framed) meromorphic
connections on $G$-bundles over the unit disc and having order two poles over
the origin (see Section \ref{sn: braids}
for full details). 
Taking the leading coefficients (irregular types) at the pole gives 
a map
$\M\to\lt_\reg$ which in fact expresses $\M$ as a fibre bundle.
Performing isomonodromic deformations of the connections then amounts
precisely to integrating
a natural flat connection on this fibre bundle (the
isomonodromy connection).
Thus, upon choosing a basepoint $A_0\in\lt_\reg$, a natural 
$P_\g$ action is obtained on the fibre $\M(A_0)$, by taking the holonomy of
the isomonodromy connection.

Now, in a previous paper \cite{smpgafm}, the author found that (for
$G=GL_n(\IC)$) the fibres $\M(A_0)$ are isomorphic to the 
Poisson Lie group $G^*$ dual to $G$ (and that the natural Poisson structures
then coincide).  
The results of Section \ref{sn: sms} enable this to be extended easily to
general $G$.
Thus isomonodromy gives a natural (Poisson) $P_\g$ action on $G^*$.

On the other hand, in their work on representations of
quantum groups at roots of unity,
De \!Concini, Kac and Procesi \cite{DKP} have written down explicitly
a Poisson action of the full braid group $B_\g=\pi_1(\lt_\reg/W)$ on $G^*$.
This was obtained by taking the classical limit of the explicit $B_\g$ 
action---the `quantum Weyl group' action---on 
the corresponding quantum group, due to
Lusztig \cite{Lus90b} and independently 
Kirillov--Reshetikhin \cite{KResh}  and Soibelman \cite{Soib}.
In this paper it will be explained how to convert the fibre bundle
$\M\to\lt_\reg$ into a bundle $\Mb\to\lt_\reg/W$ with flat connection
(and standard fibre $G^*$), by twisting by a finite group (Tits' extension
of the Weyl group by an abelian group).
Then the main result of Section \ref{sn: braids} is:

\begin{plainthm}
The holonomy action of the full braid group
$B_\g=\pi_1(\lt_\reg/W)$ on $G^*$ (obtained by integrating the flat
connection on $\Mb$) is the same as the $B_\g$ action on $G^*$ of
De \!Concini--Kac--Procesi \cite{DKP}.
\end{plainthm}

Thus the geometrical origins of the quantum Weyl group actions 
are in the geometry of meromorphic connections having order two poles. 

In Section \ref{sn: hams} a Hamiltonian description will be given of
the equations governing the isomonodromic deformations of Section 
\ref{sn: braids}. 
It will then 
be shown how this leads directly a certain flat  connection
appearing in the recent paper \cite{TL} and
featuring in the conjectural 
`Kohno--Drinfeld theorem for quantum Weyl groups'; 
see \cite{TL}, where this conjecture is explained---and 
proved for $\lsl_n(\IC)$. 
The history of this given in \cite{TL} is a little complicated:
C. De \!Concini discovered the connection and conjecture in unpublished work
around 1995. Next J. Millson and V. Toledano Laredo jointly rediscovered the
connection.
Then 
Toledano Laredo rediscovered the conjecture and found how to
prove it
for $\lsl_n(\IC)$ by translating it into the usual Kohno--Drinfeld theorem.

Our derivation of this connection of
De \!Concini--Millson--Toledano Laredo (DMT)
suggests that the theorem
of Section \ref{sn: braids} here should be interpreted as the 
classical analogue (for any $\g$) of the aforementioned conjectural
Kohno--Drinfeld theorem for quantum Weyl groups.  
The background for this interpretation comes from the paper 
\cite{Resh92} of Reshetikhin (and also \cite{Harn96, BabKit}).
In \cite{Resh92} Reshetikhin explained how
Knizhnik--Zamolodchikov type equations arise as deformations 
of the isomonodromy problem.
Although poles of order two or more are considered in \cite{Resh92},
the extra deformation parameters are not considered and so the
braiding due to the irregular types did not appear.
The derivation that will be given here of the DMT 
connection amounts
to the following statement.
If the idea of \cite{Resh92} is extended to deformations of the
isomonodromy problem for connections on $\IP^1$ with just two poles 
(of orders one and two respectively) 
then the DMT connection arises, rather than
the Knizhnik--Zamolodchikov equations.

The organisation of this paper is as follows.
Section \ref{sn: sms} swiftly states all the required results concerning the
moduli of meromorphic connections on principal $G$-bundles,
the main proofs being deferred to an appendix.
Section \ref{sn: braids} then addresses isomonodromic deformations
and proves the main theorem stated above, relating quantum Weyl group actions
to meromorphic connections.
Section \ref{sn: hams} gives the Hamiltonian approach to the 
isomonodromic deformations considered 
and shows how this leads directly to the DMT
connection.
Appendix \ref{apx: pfs} gives the proofs for Section \ref{sn: sms}.
Finally
Appendix \ref{last apx} explains how, using the results of Section 
\ref{sn: sms}, one may extend to the current setting some closely related 
theorems of a previous paper \cite{smpgafm} showing that certain
monodromy maps are Poisson.

\ 

\renewcommand{\baselinestretch}{1}		%
\small
{\em Acknowledgements.}
I would like to thank B. Dubrovin and M.S. Narasimhan for many useful 
conversations and B. Kostant for kindly supplying a proof of the fact 
that $E_8$ has no nontrivial multiplicity one representations.
The work for this paper was supported by grants from S.I.S.S.A., Trieste
and the E.D.G.E. Research Training Network HPRN-CT-2000-00101.

\renewcommand{\baselinestretch}{1.1}		%
\normalsize

\end{section}

\begin{section}{$G$-valued Stokes Multipliers} \label{sn: sms}
Let $G$ be a connected complex reductive Lie group. Fix a maximal torus
$T\subset G$ and let $\lt\subset \g$ be the corresponding Lie
algebras.
Let $\R\subset \lt^*$ be the roots of $G$ relative to $T$, so that 
as a vector space
$\g = \lt \oplus \bigoplus_{\al\in\R} \g_\al$
where $\g_\al \subset \g$ is the one-dimensional subalgebra of elements
$X\in\g$ such that $[H,X]=\al(H) X$ for all $H\in \lt$.

Let $A$ be a meromorphic connection on a principal $G$-bundle 
$P\to\Delta$ over the closed unit disc $\Delta\subset \IC$, having
a pole of order $k\ge 2$ over the origin and no others.
We view $A$ as a $\g$-valued meromorphic one-form on $P$ 
satisfying the usual conditions (\cite{KN} p.64); 
in particular the vertical component of $A$ is nonsingular. 
Upon choosing a global section $s:\Delta\to P$ of $P$ 
(which we may since every $G$-bundle over $\Delta$ is trivial), 
$A$ is determined by the $\g$-valued meromorphic one-form
$A^s:=-s^*(A)$ on $\Delta$. 
(The minus sign is introduced here simply to agree with notation 
in the differential equations literature.)
In turn $A^s=A^hdz/z^k$ for a holomorphic map $A^h:\Delta\to \g$, 
where $z$ is a fixed coordinate on $\Delta$ vanishing at
$0$.

By a {\em framing} of $P$ at $0$ we mean a point $s_0\in P_0$ of the
fibre of $P$ at $0$.
This determines the leading coefficient 
$A_0:= A^h(0)\in\g$ of $A$ independently of the choice of a
section $s$ through $s_0$.
The framed connection $(P,A,s_0)$ will be said to be 
{\em compatibly framed} if $A_0 \in \lt$. 
A compatibly framed connection is {\em generic} if 
$A_0 \in \lt_{\text{reg}}$, i.e. if $\al(A_0)\ne 0$ for all
$\al\in\R$.
Let $\G$ denote the group of holomorphic maps $g:\Delta\to G$
and let $G\flb z \frb$ be the completion at $0$ of $\G$.

\begin{lem} \label{lem: formal t}
Let $(P,A,s_0)$ be a generic compatibly framed connection with
leading coefficient $A_0$. 
Choose a trivialisation $s$ of $P$ with $s(0)=s_0$ and let 
$A^s=-s^*(A)$ as above. 
Then there is a unique formal transformation 
$\wh F\in G\flb z \frb$ and unique elements 
$A^0_{0},\ldots,A^0_{k-2},\Lambda\in \lt$ such that 
$\wh F(0) = 1$, $A^0_0=A_0$ and
$$\wh F\left[ A^0 \right] = A^s$$
where 
$A^0 := (A^0_0/z^k+\cdots+A^0_{k-2}/z^2+\Lambda/z)dz$
and 
$\wh F[A^0]$ 
denotes the gauge action 
(which, in any representation, is $\wh FA^0\wh F^{-1} + d\wh F\wh F^{-1}$).
Moreover changing the trivialisation does not change $A^0$, and
changes $\wh F$ to $\wh g\cdot\wh F$ where $\wh g\in G\flb z\frb$ is the 
Taylor expansion of some $g\in\G$ with $g(0)=1$.
\end{lem}

The proof will be given in the appendix. 
We will refer to $A^0$ as the {\em formal
type} of $(P,A,s_0)$ and to $\Lambda$ as the 
{\em exponent of formal monodromy}.
The primary aim of this section is to describe (in terms of Stokes
multipliers)
the set $\cH(A^0)$ of isomorphism classes of generic compatibly framed
connections on principal $G$-bundles over $\Delta$
with a fixed formal type $A^0$:
$$\cH(A^0) = \{  (P,A,s_0) \ \bigl\vert \ \text{formal type
$A^0$} \ \}/(\text{isomorphism}).$$

We remark that there are groups $G$ for which this description cannot be
reduced (for any $A^0$) to the $GL_n(\IC)$ case by choosing a
representation $G\subset GL_n(\IC)$ (see Lemma \ref{lem: novlem}).

Since each such principal bundle is trivial our task is equivalent to
describing the quotient 
$\{ A^s \ \bigl\vert \ \wh F[A^0]=A^s \ \text{for
some $\wh F\in G\flb z\frb$ with $\wh F(0)=1$}\}/
\{g\in \G \ \bigl\vert \ g(0)=1\}$.   
This will involve `summing' the (generally divergent) series $\wh F$
on various sectors at $0$, bounded by
`anti-Stokes directions' which are defined as follows.

Let the circle $S^1$ parameterise rays (directed lines) emanating from
$0\in\IC$. (Intrinsically this can be thought of as the boundary circle
of the real oriented blow-up of $\IC$ at $0$.)
Note that $A^0=dQ+\Lambda dz/z$ where
$Q:=\sum_{j=1}^{k-1}\frac{z^{j-k}}{j-k}A^0_{j-1}$ and
let $q:= \frac{1}{1-k}A_0z^{1-k}$ be the leading term of $Q$.
Since $A_0$ is regular, for each root $\al\in\R$, there is a non-zero 
complex number $c_\al$ such that $\al\circ q=c_\al z^{1-k}$.

\begin{defn}
The {\em anti-Stokes directions} $\IA\subset S^1$ are the directions
along which $\exp(\al\circ q)$ decays most rapidly as $z\to 0$,
i.e. the directions along which $\al\circ q(z)$ is real and negative.
\end{defn}

For $k=2$ (which will be prominent in Section \ref{sn: braids}) 
$\IA$ simply consists of the directions from $0$ to $\al(A_0)$
for all $\al\in\R$. 
(In general $\IA$ is just the inverse image under the $k-1$ fold covering
map $z\to z^{k-1}$  of the directions to the points of the set 
$\langle A_0,\R \rangle\subset \IC^*$.)
Clearly $\IA$ has $\pi/(k-1)$ rotational
symmetry and so $l:=\#\IA/(2k-2)$ is an integer.
We will refer to an $l$-tuple $\bd\subset \IA$ of 
consecutive anti-Stokes directions as a {\em half-period}.

\begin{defn} Let $d\in \IA$ be an anti-Stokes direction.

$\bullet$
The {\em roots} $\R(d)$ of $d$ are the roots $\al\in\R$
`supporting' $d$:
$$\R(d):=\{ \al\in\R \ \bigl\vert \ 
(\al\circ q)(z)\in\IR_{<0} \text{ for $z$ along } d \}.$$

$\bullet$
The {\em multiplicity} of $d$ is the number $\#\R(d)$ 
of roots supporting $d$.

$\bullet$ 
The {\em group of Stokes factors} associated to $d$ is the group
$$\ISto_d(A^0):= \prod_{\al\in \R(d)} U_\al\subset G$$
where $U_\al=\exp(\g_\al)\subset G$ is the one dimensional unipotent group
associated to $\g_\al$, and the product is taken in any order. 

$\bullet$ If $\bd\subset \IA$ is a half-period then
the {\em group of Stokes multipliers} associated to $\bd$ is
$$\ISto_\bd(A^0):= \prod_{d\in \bd} \ISto_d(A^0)\subset G.$$
\end{defn}

To understand this we note the following facts (which are proved in
the appendix):

\begin{lem}	\label{lem: s groups}
If $\bd\subset \IA$ is a half-period then 
$\bigcup_{d\in\bd} \R(d)$ is a system of positive roots in some 
(uniquely determined) root ordering.

$\bullet$ For any anti-Stokes direction $d$ the corresponding group of
Stokes factors is a unipotent subgroup of $G$ of dimension equal to the
multiplicity of $d$.

$\bullet$ For any half-period $\bd$ the corresponding group of
Stokes multipliers is the unipotent part of the Borel subgroup of $G$ 
determined by the positive roots above.

$\bullet$ The groups of Stokes multipliers corresponding to 
consecutive half-periods are the unipotent parts of opposite Borel subgroups.
\end{lem}

Now choose a sector $\Sect_0\subset \Delta$ with vertex $0$ 
bounded by two consecutive anti-Stokes directions. 
Label the anti-Stokes directions $d_1,\ldots, d_{\#\IA}$ in a positive
sense starting on the positive edge of $\Sect_0$.
Let $\Sect_i := \Sect(d_i,d_{i+1})$ denote the `$i$th sector'
(where the indices are taken modulo $\#\IA$) and define
the `$i$th supersector' to be
$\Ssect_i := \Sect\left(d_i-\frac{\pi}{2k-2},d_{i+1}+
\frac{\pi}{2k-2}\right).$
All of the sectors $\Sect_i, \Ssect_i$ are taken to be {\em open} as
subsets of $\Delta$.

\begin{thm} 	\label{thm: multisums}
Suppose $\wh F\in G\flb z \frb$ is a formal transformation as produced
by Lemma \ref{lem: formal t}. 
Then there is a unique holomorphic map 
$\Sif: \Sect_i\to G$ for each $i$ such that 

$1)$ $\Sif [A^0] =  A$, and
$2)$
$\Sif$ can be analytically continued to the supersector $\Ssect_i$ 
and then $\Sif$ is asymptotic to
$\wh F$ at $0$ within $\Ssect_i$.

Moreover,
if $t\in T$ and $g\in \G$ with $g(0)=t$ then 
$\Sigma_i( \wh g \circ \wh F \circ t^{-1})=g\circ \Sif \circ t^{-1}, $
where $\wh g$ is the Taylor expansion of $g$ at $0$.
\end{thm}
This will be proved in the appendix.
The point is that on a narrow sector there are generally many
holomorphic {isomorphisms} between $A^0$ and $A$ which are 
asymptotic to $\wh F$ and one is being chosen in a
canonical way. 

It is now easy to construct canonical $A$-horizontal 
sections of
$P$ over the sectors, using these holomorphic {isomorphisms} and
the fact that $z^\Lambda e^Q$ is horizontal for $A^0$ (which is 
viewed as a meromorphic connection on the trivial $G$-bundle).

For this we need to choose a branch of $\log(z)$ along $d_1$ which we
then  extend  in a
positive sense across $\Sect_1,d_2,\Sect_2,d_3,\ldots,\Sect_0$ 
in turn.
It will be convenient later (when $A^0$ varies) 
to encode the (discrete) choice of initial sector
$\Sect_0$ and branch of $\log(z)$ in terms of the choice of a single point 
$\wt p \in \wt\Ds$ of the universal cover of the punctured disc,
lying over $\Sect_0$.     

\begin{defn}	\label{def: sf2}
Fix data $(A^0,z,\wt p)$ as above and suppose $(P,A,s_0)$ 
is a compatibly framed connection with formal type $A^0$. 

$\bullet$
The {\em canonical fundamental solution} of $A$ on the $i$th sector is
the holomorphic map
$$\Phi_i:=\Sif z^\Lambda e^Q:\Sect_i\longrightarrow G$$
where $z^\Lambda$ uses the choice (determined by $\wt p$)
of $\log(z)$ on $\Sect_i$.

$\bullet$
The {\em Stokes factors} $K_i$ ($i=1,\ldots,\#\IA$) of 
$A$ are defined as follows.
If $\Phi_i$ is continued across the
anti-Stokes ray $d_{i+1}$ then on $\Sect_{i+1}$ we have:
$ K_{i+1} := \Phi_{i+1}^{-1}\circ\Phi_i
\text{ for } 1\le i<\#\IA 
\text{ and }
 K_{1} := \Phi_{1}^{-1}\circ\Phi_{\#\IA} \circ M^{-1}_0,$
where $M_0:=e^{2\pi\sqrt{-1}\cdot\Lambda}\in T$ is the `formal monodromy'.

$\bullet$
The {\em Stokes multipliers} $S_i$ ($i=1,\ldots,2k-2$) of 
$A$ are 
\begin{equation}	\label{eqn: sm}
S_i:= K_{il}\cdots K_{(i-1)l+1},
\end{equation}
where $l=\#\IA/(2k-2)$. Equivalently if 
$\Phi_{il}$ is continued across 
$d_{il+1},\ldots,d_{(i+1)l}$ and 
onto $\Sect_{(i+1)l}$ then:
$\Phi_{il}=\Phi_{(i+1)l}S_{i+1}$
for $i=1,\ldots,2k-3$, and
$\Phi_{il}=\Phi_lS_1M_0$ 
for $i=2k-2$.
\end{defn}
Note that the canonical solutions $\Phi_i$ are appropriately
equivariant under change of trivialisation, so are naturally
identified with $A$-horizontal sections of $P$. 
It follows that the Stokes factors and Stokes multipliers 
are constant ($z$-independent) elements  of $G$.
Also, from
the proof of Lemma \ref{lem: s groups}, note that $S_i$ 
uniquely determines each 
Stokes factor appearing in \eqref{eqn: sm}.
In the appendix we will establish the basic lemma:
\begin{lem}	\label{lem: sf in sfg}
$K_i\in \ISto_{d_i}(A^0)$ and $S_j\in\ISto_{\bd}(A^0)$
where $\bd= (d_{(j-1)l+1},\ldots,d_{jl})$.
\end{lem}
It is immediate from Definition \ref{def: sf2} and 
the last part of Theorem \ref{thm: multisums}
that the Stokes multipliers are
independent of the trivialisation choice in Lemma \ref{lem: formal t},
and so are well defined (group valued) functions on $\cH(A^0)$.
The main result of this section is then:

\begin{thm}   \label{thm: irr RH im}
Fix the data $(A^0,z,\wt p)$ as above.
Let $U_+=\ISto_\bd(A^0)$ where $\bd=(d_1,\ldots,d_l)$ is the first
half-period and let $U_-$ denote the opposite full unipotent subgroup
of $G$. 
Then the `irregular Riemann--Hilbert map' taking the Stokes multipliers 
induces a
{\em bijection}
$$\cH(A^0)\mapright{\cong} (U_+\times U_-)^{k-1}; \qquad 
[(P,A,s_0)]\longmapsto (S_1,\ldots,S_{2k-2}).$$
In particular $\cH(A^0)$ is isomorphic to a complex vector space 
of dimension $(k-1)\cdot(\#\R)$.
\end{thm}
\pf
For injectivity, suppose we have two compatibly framed 
meromorphic connections with 
$\wh F_1[A^0]=A^s_1$ and $\wh F_2[A^0]=A^s_2$ and
having the same Stokes multipliers.
Therefore the Stokes factors are also equal and it follows immediately that
$\Sigma_i(\wh F_2)\circ\Sigma_i(\wh F_1)^{-1}$ has no monodromy
around $0$ and does not depend on $i$,
and thereby defines a holomorphic map $g:\Ds\to G$.
Thus on any sector $g$ has asymptotic expansion $\wh F_2\circ \wh F_1^{-1}$
and so (by Riemann's removable singularity theorem) we deduce 
the formal series $\wh F_2\circ \wh F_1^{-1}$ is actually convergent
with the function $g$ as sum. 
This gives an isomorphism between the connections we began
with: they represent the same point in $\cH(A^0)$.
Surjectivity follows from the $G$-valued analogue of a theorem of Sibuya,
which we will prove in the appendix.
\epf

To end this section we will show that 
the Stokes multipliers of a holomorphic family of
connections vary holomorphically with the parameters of the family.
Suppose we have a family of compatibly framed meromorphic connections on
principal $G$-bundles over the disc $\Delta$, 
parameterised by some polydisc $X$.
Upon choosing compatible trivialisations this family may be written as
$$ A^s = A^h \frac{dz}{z^k}$$
for a holomorphic map $A^h:\Delta\times X \to \g$ with leading coefficient 
a holomorphic map $A_0=A^h\bigl\vert_{z=0}:X\to \lt_\reg$.
The proof of Lemma \ref{lem: formal t} is completely algebraic and remains
unchanged upon replacing the coefficient ring $\IC$ by the ring 
$\cO(X)$ of holomorphic functions on $X$; there
is a unique formal transformation 
$\wh F\in G\bigl(\cO(X)\flb z \frb\bigr)$ and unique holomorphic maps 
$A^0_{0},\ldots,A^0_{k-2},\Lambda : X \to \lt$ such that 
$\wh F\bigl\vert_{z=0} = 1$, $A^0_0=A_0$ and
$\wh F\left[ A^0 \right] = A^s$
where 
$A^0 := (A^0_0/z^k+\cdots+A^0_{k-2}/z^2+\Lambda/z)dz$.
Given a point $x\in X$ let $\wh F_x\in G\flb z\frb$ denote the 
corresponding formal bundle automorphism.

Now choose any basepoint $x_0\in X$ and let $\IA_0\subset S^1$ 
denote the anti-Stokes
directions associated to $A_0(x_0)$.
Let $\check S\subset \Delta$ 
be any  sector (with vertex $0$) and
whose closure contains none of the directions in $\IA_0$. 
By continuity there is a neighbourhood $U\subset X$ of $x_0$ such that
none of the anti-Stokes directions associated to $A_0(x)$ lie in $\check S$, 
for any $x\in U$.
We will always label the sectors such that $\check S\subset \Sect_0$.

\begin{lem} \label{lem: holom dpce}
In the situation above
the holomorphic maps $\Si_0(\wh F_x) : \check S\to G$ (defined for each $x$ 
in Theorem \ref{thm: multisums} and restricted to $\check S$) 
vary holomorphically with $x\in U$ and so constitute a holomorphic map
$$\Si_0(\wh F):\check S\times U\to G.$$
\end{lem}

This will be proved in the appendix. 

\begin{cor} \label{cor: hol dpce}
In the situation above, 
taking Stokes multipliers defines a holomorphic map 
$U\to (U_+\times U_-)^{k-1}$ from the parameter space $U$ to the space of 
Stokes multipliers.
In particular if $A^0$ is any formal type then $\cH(A^0)$ is a coarse 
moduli space in the analytic category.
\end{cor}
\pf Lemma \ref{lem: holom dpce} implies each of the sums 
$\Si_{il}(\wh F_x)$ varies holomorphically with $x\in U$
(even though the integer $l$ may jump; $\Si_{il}(\wh F_x)$ is defined 
invariantly as the `sum' of $\wh F_x$ on  the sector
$\check S\cdot\exp\bigl(\frac{i\pi\sqrt{-1}}{k-1}\bigr)$).
Thus, once a branch of $\log(z)$ is chosen on $\check S$, 
the canonical solutions 
$\Phi_{il}$ also vary holomorphically with $x$.
The Stokes 
multipliers are defined directly in terms of these canonical solutions  and so
also vary holomorphically.
That $\cH(A^0)$ is a coarse moduli space is immediate from this and 
Theorem \ref{thm: irr RH im}.
\epf

\end{section}

\begin{section}{Isomonodromic Deformations}  \label{sn: braids}

In this section we will define and study isomonodromic deformations 
of generic compatibly framed meromorphic connections on principal $G$-bundles
over the unit disc, having an order two pole at the origin.
Due to the results of the previous section the definition
 is now a straightforward matter. 
(The $GL_n(\IC)$ case over $\IP^1$, with arbitrary many poles of arbitrary
order was defined in \cite{JMU81} and studied further in \cite{smid}.)

The main aim here is to describe a relationship between 
isomonodromic deformations and certain braid group actions
arising in the theory of quantum groups.
In brief this relationship is as follows.
In \cite{smpgafm} the author identified the Poisson Lie group $G^*$ dual to
$G=GL_n(\IC)$ with a certain moduli space $\M(A_0)$ of meromorphic
connections on vector bundles (principal $GL_n(\IC)$ bundles) 
over the unit disc and having an order two pole
at the origin and `irregular type' $A_0\in\lt_\reg$. 
The previous section enables us to extend this identification 
easily to arbitrary $G$.

By considering `isomonodromic deformations' of such connections 
(where $A_0$ plays the role of deformation parameter) one obtains
an action of the pure braid group $P_\g=\pi_1(\lt_\reg)$ on 
$\M(A_0)\cong G^*$.
This is purely geometrical: there is a moduli space $\M$ of meromorphic
connections fibring over $\lt_\reg$ (with fibre $\M(A_0)\cong G^*$ over
$A_0$) and having a natural flat (Ehresmann)
connection---the isomonodromy connection. The $P_\g$ action is just the
holonomy of this connection.

On the other hand 
De \!Concini--Kac--Procesi \cite{DKP} have described explicitly
an action of the full braid group $B_\g=\pi_1(\lt_\reg/W)$ on $G^*$ in their
work on representations of quantum groups at roots of unity.
(This is for simple $\g$, which is certainly the most interesting case.)
This action is the classical version of the 
`quantum Weyl group'
actions  of $B_\g$ on a quantum group (the quantisation of $G^*$) which were
defined by
Lusztig, Kirillov--Reshetikhin, and Soibelman (see \cite{DP1565} for more
details; in particular Section 12 gives the definition of the quantum
group having classical limit $G^*$).

In this section we will explain how to convert $\M$ into a fibre bundle
$\Mb\to\lt_\reg/W$ with flat connection, using Tits' `extended Weyl group'
\cite{Tits66}, 
and then prove that the holonomy action of $B_\g$ on the fibres of $\Mb$
(which are still isomorphic to $G^*$) 
is precisely the action of De \!Concini--Kac--Procesi.
Thus we have a geometrical description of their action;
roughly speaking 
the infinite part (related to $P_\g$) 
of the $B_\g$ action comes 
from geometry whereas the rest (related to the Weyl group) 
is put in by hand.

We have restricted to the order two pole case over a disc here 
since that is what 
is required for the application we have in mind here.
However the results of Section \ref{sn: sms} do immediately facilitate the
definition of isomonodromic deformations in much more generality. 

\ 

{\bf The fibration $\M\to\lt_\reg$.}\ \ \ 
Fix a connected complex simple Lie group $G$ and a maximal torus
$T\subset G$. 
In terms of the definitions of Section \ref{sn: sms} we have:

\begin{defn}
The moduli space $\M$ is the set of isomorphism classes of triples
$(P,A,s_0)$ of generic compatibly framed meromorphic connections $A$
on principal
$G$-bundles $P\to\disc$ having an order two pole at the origin.
\end{defn}

Denote by $\pi:\M\to \lt_\reg$ the surjective map taking the leading
coefficient 
of the connections in $\M$.
Let $\M(A_0)\subset \M$ be the
fibre of $\pi$ over the point $A_0\in \lt_\reg$.
($A_0$ will be called the `irregular type'; it determines
the irregular part of the formal type $A^0$ of a
connection in $\M$.) 

\begin{prop} \label{prop: imd conn}
$\M$ has the structure of complex analytic fibre bundle over $\lt_\reg$ with 
standard fibre $U_+\times U_-\times \lt$, where $U_\pm$ are the unipotent
parts of a pair of opposite Borel subgroups $B_\pm\subset G$ containing $T$.

Moreover there is a canonically defined flat (Ehresmann) connection on 
$\M\to\lt_\reg$; the isomonodromy connection.
\end{prop}

\pf
Fix an irregular type $A_0\in\lt_\reg$. 
This determines anti-Stokes directions at $0$ as in Section \ref{sn: sms}.
Choose $\wt p\in \wt\Ds$ as in Definition \ref{def: sf2}, 
determining an initial sector $\Sect_0$ and branch of $\log(z)$.
Then if we define $U_\pm$ in terms of the first half-period as in
Theorem \ref{thm: irr RH im}, this choice determines an isomorphism
\begin{equation}	\label{eqn: mao im}
\M(A_0)\cong U_+\times U_-\times \lt
\end{equation}
as follows.
There is a surjective map $\M(A_0)\to\lt$ taking a connection to its exponent
of formal monodromy $\Lambda$ (the residue of its formal type), 
as defined in Lemma \ref{lem: formal t}.
By definition the fibre of this map over $\Lambda\in\lt$ is 
$\cH(A^0)$ where $A^0:= \left(A_0/z^2+\Lambda/z\right)dz$.
Then by Theorem \ref{thm: irr RH im} each such fibre is canonically isomorphic
to $U_+\times U_-$ (using the choice of $\wt p$ made above) and so 
\eqref{eqn: mao im} follows.

Now if we vary $A_0$ slightly, since the anti-Stokes directions depend
continuously on $A_0$ and $\Sect_0$ is open, 
we may use the same $\wt p$ for all $A_0$ in some 
neighbourhood of the original one.
The above procedure then gives a local trivialisation of $\M\to\lt_\reg$ 
over this neighbourhood, implying it is indeed a fibre bundle.

If we repeat this for each $A_0\in\lt_\reg$ and each choice of $\wt p$ we
obtain an open cover of $\lt_\reg$ with a preferred 
trivialisation of $\M$ over each open set.
The clutching maps for this open cover are clearly constant 
(involving just rearranging the Stokes factors into Stokes multipliers 
in different ways
and conjugating by various exponentials of $\Lambda$),
and so we have specified
a flat connection on the fibre bundle $\M\to\lt_\reg$, the local horizontal
leaves of which contain meromorphic 
connections with the same Stokes multipliers and exponent
of formal monodromy (for some---and thus any---choice of $\wt p$).
\epf

\begin{rmk}
The isomonodromy connection may be viewed profitably as an analogue of the 
Gauss--Manin connection in non-Abelian cohomology (which has been 
studied by Simpson \cite{Sim94ab}). 
Extending Simpson's terminology we will call the 
above definition the `Betti' approach to isomonodromy. 
There is also an equivalent `De \!Rham' approach 
involving flat meromorphic connections on $G$-bundles over products 
$\Delta\times U$ for open neighbourhoods $U\subset\lt_\reg$.
(This is well-known to isomonodromy experts in the $GL_n(\IC)$ case.)
This point of view has been described by the author in 
\cite{smid} Section 7 for the $GL_n(\IC)$ case; this now 
extends immediately to arbitrary $G$ (see \cite{smid} Theorem 7.2 
in particular for the De \!Rham approach).
\end{rmk}

The next step is to convert $\M$ into a fibre bundle 
$\Mb\to\lt_\reg/W$ with flat connection (where
$W:= N(T)/T$ is the Weyl group), so that one obtains a
holonomy action of the full braid group 
$B_\g:= \pi_1(\lt_\reg/W)$ on the fibres, rather than just an action of the
pure braid group $P_\g:= \pi_1(\lt_\reg)$.
(This step is closely related to a similar step taken by
Toledano Laredo in \cite{TL}.)
One would like simply to quotient $\M$ by an action of $W$ covering the
standard free action on $\lt_\reg$. 
Indeed if there was a homomorphic section $W\to N(T)\subset G$ of the
canonical projection 
$\pi_N: N(T)\to W$, then we could simply act on $\M$ by constant
gauge transformations. 
However there is no such section in general, even for $SL_2(\IC)$.
(For $GL_n(\IC)$ one may use the section given by `permutation matrices' but
here we require a general approach.)

The standard way around this problem was found by Tits \cite{Tits66};
there is a finite abelian extension 
$$1\to \Ga_1 \to \Ga \mapright{\pi_\Ga} W \to 1$$ 
of $W$ (where $\Ga$ is finite and $\Ga_1$ is abelian) and an inclusion 
$\iota: \Ga \hookrightarrow N(T)$ covering the identity in $W$ (i.e. so that 
$\pi_N\circ \iota = \pi_\Ga$). The group $\Ga$ is Tits' extended Weyl group.

\begin{rmk} \label{rmk: ewg}
A construction of $\Ga$ is as follows.
Choose a positive Weyl chamber,
label the simple roots by $i=1,\ldots,n$ and choose Chevalley generators
$\{e_i, f_i, h_i\}$ of $\g$ as usual.
Let
$$t_i := \exp(f_i)\exp(-e_i)\exp(f_i)\in G.$$
One then knows (from \cite{Tits66}) that 1) 
these $t_i$ satisfy the braid relations for $\g$ and so
determine a homomorphism 
$B_\g\to G$, and 2) the image $\Ga$ of $B_\g$ in $G$ has the
properties stated above. 
(We note for later use that replacing $t_i$ by $t_i^{-1}$ here determines
another homomorphism $B_\g\to G$ with the same image.)
\end{rmk}

We could now act with $\Ga$ on $\M$ by gauge transformations, but then the
quotient
would not be a fibre bundle over $\lt_\reg/W$, since this action is not free 
(e.g. $\Ga_1$ acts trivially on formal types, but non-trivially on other 
connections).
To get around this we first pull back $\M\to \lt_\reg$ to the 
Galois $\Ga_1$ cover $\wh\lt_\reg$ of $\lt_\reg$.
(In other words $\wh\lt_\reg:= \wt\lt_\reg/K$, where 
$\wt\lt_\reg$ is the universal cover of $\lt_\reg$ 
and $K:=\ker(B_\g\to \Ga)=\ker(P_\g\to \Ga_1)$.)   
Then define $\wh \M := \pr^*(\M)$ to be the pullback of the bundle $\M$ along
the covering map 
$\pr:\wh\lt_\reg\to\lt_\reg$. The connection on $\M$ pulls back
to a flat connection on $\wh \M\to\wh\lt_\reg$.

Finally we can now act with $\Ga$ on $\wh \M$ by gauge transformations,
covering the canonical free action of $\Ga$ on $\wh\lt_\reg$, to obtain a
fibre bundle 
$\Mb := \wh \M/\Ga \to \lt_\reg/W$.
In summary we have the commutative diagram:
\begin{equation}	\label{cd: Ms}
\begin{array}{ccccc}
  \wh\M & \mapright{} & \M & & \Mb \\
\mapdown{} & & \mapdown{} & & \mapdown{}  \\
  \wh\lt_\reg & \mapright{\pr} & \lt_\reg & \mapright{} & \lt_\reg/W, 
\end{array}
\end{equation}
where the horizontal maps are finite covering maps and the vertical maps are
fibrations. 

\begin{lem} \label{lem: conn descends}
The connection on $\wh\M=\pr^*(\M)$ is $\Ga$ invariant and so descends 
to a flat connection on $\Mb\to\lt_\reg/W$.
\end{lem}
\pf
Choose $g\in\Ga\subset N(T)$ and $\wh A_0\in\wh\lt_\reg$.
It is sufficient to show that, under the action of $g$,
local horizontal sections of
$\wh \M$ over a neighbourhood of $\wh A_0$ become 
horizontal sections 
over a neighbourhood of $g(\wh A_0)$.
To this end, choose open $U\subset\lt_\reg$ containing $A_0:=\pr(\wh A_0)$ 
as in the
proof of Proposition \ref{prop: imd conn} and so small that $\pr^{-1}(U)$
consists of $\#\Ga_1$ connected components. 
Let $\wh U$ be the component containing $\wh A_0$.

Now choose $\wt p\in \wt\Ds$ as in Proposition \ref{prop: imd conn} and thereby
obtain a (horizontal) trivialisation 
of $\M$ over $U$ and also of $\wh\M$ over $\wh U$:
$$\wh\M\bigl\vert_{\wh U} \cong U_+\times U_-\times\lt\times \wh U,$$
where $U_\pm$ are determined by $A_0$ and $\wt p$ as in 
Theorem \ref{thm: irr RH im}.
Let $A'_0 := gA_0g^{-1}$ so that $\pr(g(\wh A_0)) = A'_0$. 
Thus $\pr(g(\wh U))$ is a neighbourhood of $A'_0$ over which we may
trivialise $\M$
using the same choice of $\wt p$ as above (since $A_0$ and $A'_0$ determine
the same set of anti-Stokes directions). 
Thus in turn 
$\wh\M\bigl\vert_{g(\wh U)} \cong U'_+\times U'_-\times\lt\times g(\wh U),$
where $U'_\pm$ are determined by $A'_0$ and $\wt p$ as in 
Theorem \ref{thm: irr RH im}.

Finally we claim that $U'_\pm=gU_\pm g^{-1}$ and that, 
in terms of the above trivialisations, the action of
$g$ on $\wh \M$ is given by 
\begin{equation}	\label{eqn: gamma action}
g(S_+,S_-,\Lambda) = (gS_+g^{-1}, gS_-g^{-1}, g\Lambda g^{-1})
\end{equation}
(together with the standard action on the base $\wh \lt_\reg$),
where $(S_+,S_-,\Lambda)\in U_+\times U_-\times\lt$. 
Since there is no dependence on the base, 
this clearly implies the proposition.
The claim is established by a straightforward unwinding of the definitions. 
\epf

In \cite{smpgafm} (for $G=GL_n(\IC)$)  
it was found to be natural to identify the space
$U_+\times U_-\times \lt$ of monodromy data with the simply-connected Poisson
Lie group $G^*$ dual to $G$, which we will now do here in general
(cf. also Appendix \ref{last apx} for motivation).
Given a choice $B_\pm$ of opposite Borel subgroups of $G$ with
$B_+\cap B_-=T$, the group $G^*$ is defined to be
\begin{equation}	\label{gstar}
 G^* := \{ (b_-,b_+,\Lambda)\in B_-\times B_+\times\lt \ \bigl\vert \  
\delta_-(b_-)\delta_+(b_+)=1, \delta_+(b_+)=\exp(\pi i \Lambda) \}, 
\end{equation}
where $\delta_\pm:B_\pm\to T$ is the natural projection (with kernel the
unipotent part $U_\pm$ of $B_\pm$) and $\exp:\lt\to T$ is 
the exponential map for $T$.
This is a simply-connected
(indeed contractible) subgroup of $B_-\times B_+\times\lt$ (where
$\lt$ is a group under $+$) of the same dimension as $G$.

The group $G^*$ is then identified with $U_+\times U_-\times \lt$ 
as follows (cf. \cite{smpgafm} Definition 20)
\begin{equation} \label{eqn: im to gs}
U_+\times U_-\times \lt\cong G^*;\qquad 
(S_+,S_-,\Lambda)\mapsto (b_-,b_+,\Lambda)
\end{equation}
where $b_- = e^{-\pi i \Lambda} S_-^{-1}$ and 
$b_+ = e^{-\pi i \Lambda} S_+e^{2\pi i\Lambda}$, 
so that $b_-^{-1}b_+ = S_-S_+\exp(2\pi i\Lambda)$.

Thus the fibrations \eqref{cd: Ms} can now be viewed as having standard fibre 
$G^*$ (although they are not principal $G^*$-bundles).

The final (trivial) complication is that we have $G^*$ simply-connected,
whereas \cite{DKP, DP1565} use the quotient group
defined by omitting
the $\Lambda$ component in \eqref{gstar} (or equivalently 
one only remembers $e^{\pi i \Lambda}$ rather than $\Lambda$). 
We will abuse notation and denote both groups $G^*$; in terms of the braid
groups actions this is reasonable since 
1) It is immediate that the connection on $\Mb$ is invariant
under the corresponding action of the lattice $\ker(\exp(\pi i\cdot):\lt\to G)$
so descends to give a flat connection on the quotient bundle 
(still denoted $\Mb$), and
2) The $B_\g$ action of \cite{DKP} lifts to an action on our $G^*$
simply by acting on $\Lambda$ via the standard Weyl group action.

The main result is then:

\begin{thm} \label{thm: holonomy}
The holonomy action of the full braid group
$B_\g=\pi_1(\lt_\reg/W)$ on $G^*$ (obtained by integrating the flat
connection on $\Mb$) is the same as the $B_\g$ action on $G^*$ of
De \!Concini--Kac--Procesi \cite{DKP}.
\end{thm}
\pf
Choose a real basepoint $A^*_0\in\lt_{\IR,\reg}\subset\lt_\reg$ in the
Weyl chamber chosen in Remark \ref{rmk: ewg} above 
(cf. also \eqref{eq: weyl chambers}).
Then the corresponding set $\IA$ 
of anti-Stokes directions consists of just the two
halves of the real axis. 
Let $\Sect_0$ be the lower half disc, choose a point $p$ on the 
negative imaginary axis and let $\wt p\in \wt\Ds$ 
be the point lying over $p$ and on the branch of 
logarithm having $\log(-i)=3\pi i/2$. 
Define the positive roots $\R_+$, the groups $B_\pm,U_\pm$ and in turn $G^*$
to be those determined by $A^*_0$ and $\wt p$.
These choices determine an isomorphism $\M(A^*_0)\cong G^*$ via Theorem 
\ref{thm: irr RH im} and \eqref{eqn: im to gs}.
(One may check $\R_+$ is the set of positive roots corresponding to the
chosen positive Weyl chamber.)

Now for each simple root $\al=\al_i\in\R_+$ Brieskorn \cite{Bries71}
defines the following path 
$\ga_i$ in $\lt_\reg$.
Let $s_i$ be the complex reflection acting on $\lt$ corresponding to
$\al$ (the reflection fixing the hyperplane $\ker(\al)$ and respecting 
the Killing form).
Let $L_i$ be the complex line in $\lt$ containing $A^*_0$ and 
$A'_0:=s_i(A^*_0)$ and let $I_i$ be the real line segment from
$A^*_0$ to $A'_0$.
Then define the path $\ga_i:[0,1]\to L_i$ from $A^*_0$ to
$A'_0$ such that $[0,\frac{1}{3}]\cup [\frac{2}{3},1]$ 
maps to $I_i$ and $[\frac{1}{3},\frac{2}{3}]$ 
maps to a small semi-circle turning in a positive sense and centred on the
midpoint of $I_i$.
According to \cite{Bries71}, if the semi-circles are sufficiently small,
these paths $\ga_i$ are in $\lt_\reg$ and descend to loops in $\lt_\reg/W$ 
representing generators of $\pi_1(\lt_\reg/W)$.

For our purposes here we choose the above semi-circles
so small that, as $A_0$ moves along $\ga_i$, precisely one anti-Stokes 
direction crosses over the point $p$---an anti-Stokes direction supported
just by $-\al$ and moving in a positive sense.
(To see this is possible observe that for any $\be\in\R_+\setminus\{\al\}$,
$\be(A^*_0)$ and $\be(A'_0)$ are real and positive, since 
$\al$ is the only positive root made negative by $s_i$.
Thus by linearity $\be(I_i)\subset\IR_{>0}$.
Hence if $\ga_i$'s
semi-circle is sufficiently small $\be(A_0)$ does not cross the imaginary axis
for any $A_0$ on $\ga_i$ and therefore no anti-Stokes direction supported by 
$\pm\be$ crosses $p$. 
Finally observe that, as $A_0$ moves along $\ga_i$,
$\al(A_0)$ starts in $\IR_+$, moves towards $0$, makes a positive semi-circle
around $0$ then moves away from $0$ along  $\IR_-$. Since $q=-A_0/z$ here 
this implies the anti-Stokes direction supported by $-\al$ crosses $p$,
and the one supported by $\al$ crosses the positive imaginary axis.)

We now wish to calculate the holonomy isomorphism 
$\M(A^*_0)\cong \M(A'_0)$ obtained by integrating the isomonodromy
connection along the path $\ga_i$.
As in Proposition \ref{prop: imd conn} we have canonical 
descriptions of the fibre of $\M$ over both $A^*_0$ and $A'_0$:
$$\M(A^*_0)\cong U_+\times U_-\times \lt\qquad
\M(A'_0)\cong U'_+\times U'_-\times \lt$$
using the chosen $\wt p$ in both cases, where
(as in Lemma \ref{lem: conn descends}) $U'_\pm=gU_\pm g^{-1}$ 
for any $g\in\Ga$ with $\pi_\Ga(g)=s_i\in W$.
Thus we want to find the corresponding 
isomorphism $U_+\times U_-\times \lt$ $\cong$ $U'_+\times U'_-\times \lt$.
To describe it we will need the following maps.
Let $U_i=\exp(\g_{\al})$ be the root group corresponding to the simple
root $\al=\al_i$.
Then there is a homomorphism $$\xi_i:U_+\mapright{} U_i$$
with the property that
if any $S\in U_+$ is factorised (in any order) as a product of elements
$u_\be\in U_\be$ for $\be\in\R_+$ (with each $\be$ appearing just once)
then $u_{\al}=\xi_i(S)$.
(The existence of $\xi_i$ may be seen as follows:
The set $\Psi:=\R_+\setminus\{\al\}$ is a closed set of roots so
$U_\Psi:=\prod_{\be\in\Psi}U_\be$ is a subgroup of $U_+$.
By \cite{Bor91} Proposition 14.5(3) $U_\Psi$ is a normal subgroup. 
(It is sufficient to prove $U_i$ normalises $U_\Psi$.)
Then $\xi_i$ is taken to be the projection $U_+\to U_+/U_\Psi$ where
$U_i\cong U_+/U_\Psi$ via the inclusion $U_i\subset U_+$.)
Similarly we have maps $\xi_{-i}:U_-\to U_{-\al}$. 

\begin{prop} \label{prop: hol im}
The holonomy isomorphism $\M(A^*_0)\cong\M(A'_0)$ induced by the
isomonodromy connection is given by
$U_+\times U_-\times\lt\to U'_+\times U'_-\times\lt;
(S_+,S_-,\Lambda)\mapsto(S'_+,S'_-,\Lambda)$
where 
\begin{equation}	\label{eqn: hol im in SMs}
S'_+ := \xi_i(S_+)^{-1}S_+M_0\xi_{-i}(S_-)M_0^{-1},\qquad
S'_- := \xi_{-i}(S_-)^{-1}S_-\xi_i(S_+),
\end{equation}
and $M_0:=\exp(2\pi i \Lambda)$.
\end{prop}
\pf
We must find the transition maps between the local trivialisations just before
and after the anti-Stokes direction $d_{-\al}$ supported by $-\al$ crosses $p$.
By perturbing $\wt p$ (and therefore also $p$) slightly this is equivalent to
finding the transition map relating the two situations appearing in 
Figure \ref{bafi fig}, where $p$ moves but $A_0$---and thus all the anti-Stokes
directions---remain fixed.
That is, we must find the composite map 
\begin{equation} \label{eqn: composite}
U_+\times U_-\times\lt\overset{}{\cong}
\M(A_0)\overset{}{\cong} U'_+\times U'_-\times\lt
\end{equation}
where the first (resp. second) isomorphism is determined by the $\wt p$ choice
in the left (resp. right) diagram in Figure \ref{bafi fig}.

\setcounter{figure}{0}
\begin{figure}[hbt] 	
\centerline{\input{bafi.pstex_t}}	
\caption{} \label{bafi fig}
\end{figure}

Choose arbitrary $(S_+,S_-,\Lambda)\in U_+\times U_-\times\lt$ and let 
$A$ be a connection on the trivial $G$-bundle over $\Delta$ with
isomorphism class in $\M(A_0)$ corresponding to 
$(S_+,S_-,\Lambda)$ under the left-hand isomorphism in \eqref{eqn: composite}.
Let $\Phi_0,\ldots, \Phi_3$ and $\Psi_0,\ldots \Psi_3$ be the canonical 
fundamental solutions of $A$ on the sectors indicated in Figure \ref{bafi fig}.
(Except for $\Phi_0$ the indexing of these differs from Definition 
\ref{def: sf2}.)
Since the $\log(z)$ choice on the sector containing $p$ is extended to the
other sectors in a negative sense, we immediately deduce
$\Phi_3=\Psi_3, \Phi_2=\Psi_2, \Phi_1= \Psi_1$ and
$\Phi_0=\Psi_0 M_0$.

Now let $K_\pm$ denote the Stokes factors of $A$ across $d_{\pm\al}$ in
the left diagram in Figure \ref{bafi fig}
(and similarly $K'_\pm$ for the right diagram). 
Clearly $K_+=K'_+$ since across $d_\al$ we have
$$K_+:= \Phi_2^{-1}\Phi_1 = \Psi_2^{-1}\Psi_1 =: K'_+.$$
Across $d_{-\al}$,  $K_-:= \Phi_0^{-1}\Phi_3$ and $K'_-:= \Psi_0^{-1}\Psi_3$ 
so that $K'_-=M_0K_-M_0^{-1}$.

In turn the Stokes multipliers are defined by the equations
\begin{align*}
\Phi_2&=\Phi_0S_- & \Phi_0&=\Phi_2S_+ \\
\Psi_1&=\Psi_3S'_- & \Psi_3&=\Phi_1S'_+
\end{align*}
where in the left/right column the fundamental solutions are continued into
the 
left/right half-plane before being compared, respectively.
(Here we prefer to index the Stokes multipliers by $+$ and $-$ rather
than $1$ and $2$ as in Section \ref{sn: sms}.)
Combining this with the above expression for the Stokes factors we deduce
$$
S'_+=K^{-1}_+S_+M_0K_-M_0^{-1},\qquad
S'_-=K_-^{-1}S_-K_+.$$
Finally, from the alternative definition of the Stokes multipliers in terms of
Stokes factors in Definition \ref{def: sf2}, we find
$\xi_{\pm i}(S_\pm)=K_\pm$ thereby completing the proof of the proposition.
\epf

To rewrite this holonomy isomorphism in terms of the Poisson Lie groups it is 
convenient to introduce the following notation.
If $b_\pm = v_\pm t^{\pm 1}= t^{\pm 1}u_\pm$ where $t\in T$ and 
$u_\pm,v_\pm\in U_\pm$ then
$$\pibpm := \xi_{\pm i}(v_\pm)^{-1}, \qquad 
b_\pm^i := \xi_{\pm i}(u_\pm)^{-1}.$$
(The inverted left and right $\pm\al$ components of $b_\pm$ respectively.)
Under the identification \eqref{eqn: im to gs}, the isomorphism 
\eqref{eqn: hol im in SMs} then simplifies to
\begin{equation} \label{eqn: hol im in gstar1}
(b_-,b_+,\Lambda)\,\mapsto\, 
(\pibp b_- b_-^i,\, \pibp b_+ b_-^i,\, \Lambda).
\end{equation}

Clearly we may quotient by the lattice $\ker(\exp(\pi i\cdot):\lt\to G)$
(i.e. forget the $\Lambda$ component above) since $\Lambda$ only appears as 
$e^{2\pi i\Lambda}$ in the formulae and 
$t:=e^{\pi i \Lambda}=\delta_+(b_+)$ is retained.

Now if we choose $\wh A_0^*\in \pr^{-1}(A_0^*)$ and lift $\ga_i$ canonically
to a path $\wh\ga_i$ in $\wh\lt_\reg$ starting at $\wh A_0^*$,
then the holonomy of the connection on $\wh \M$ along $\wh\ga_i$ is also
given by Proposition \ref{prop: hol im} (since the connection is pulled back
from $\M$). 
Then quotienting by $\Ga$ enables us to identify the fibres
$\wh\M(\wh A_0^*)$ and $\wh\M(\wh\ga_i(1))$ via the gauge action of $t_i$.
(This uses the fact that the element of $B_\g$ determined by $\ga_i$ maps to
$t_i$ under the surjection $B_\g\to\Ga$.) The $\Ga$ action on Stokes 
multipliers
was given in \eqref{eqn: gamma action}, and so we deduce the following 
formula 
for the holonomy isomorphism $G^*\to G^*$ for the connection 
on $\Mb$ around the loop $\ga_i/W$:
\begin{equation} \label{eqn: hol im in gstar2}
(b_-,b_+)\,\mapsto\, 
(t_i^{-1}\pibp b_- b_-^it_i,\, t_i^{-1}\pibp b_+ b_-^it_i).
\end{equation}

Finally we must compare \eqref{eqn: hol im in gstar2} with 
the generators of the braid group action of De \!Concini--Kac--Procesi 
\cite{DKP}.
In \cite{DKP} the braid group $B_\g$ is defined abstractly by generators and
relations, rather than as a fundamental group.
Namely one has generators $T_i$ (one for each simple root $\al_i$) and
relations 
$$T_iT_jT_i\cdots = T_jT_iT_j\cdots$$
for $i\ne j$, 
where the number of factors on each side equals the order of the element
$s_is_j$ of the Weyl group.
The action of $B_\g$ on $G^*$ is given in Section 7.5 of \cite{DKP} 
by the following formula: 
$$T_i(t^{-1}u_-^{-1},tu_+)\ =\ 
\bigl(t_it^{-1}(u_-^{(i)})^{-1}(\exp \tilde x_ie_i)t_i^{-1},\, 
t_it^{-1}(\exp \tilde y_i f_i)t^2u_+^{(i)}t_i^{-1}\bigr)$$
where, in our notation, 
$\exp \tilde x_ie_i = \xi_i(u_+)^{-1}$, 
$\exp \tilde y_if_i = \xi_{-i}(u_-)$, 
$u_+^{(i)}=u_+\exp \tilde x_ie_i$, and
$u_-^{(i)}=u_-(\exp \tilde y_if_i)^{-1}$.
One may readily check this is the same as 
\begin{equation} \label{eqn: hol im in gstar3}
(b_-,b_+)\,\mapsto\, 
(t_i \pibm b_- b_+^it_i^{-1},\, t_i\pibm b_+ b_+^it_i^{-1})
\end{equation}
where $(b_-,b_+)=(t^{-1}u_-^{-1},tu_+)\in G^*$.
In turn it is straightforward to check this is precisely the inverse map
to \eqref{eqn: hol im in gstar2}.
Thus, if we choose to identify the (abstractly presented) braid group with
$\pi_1(\lt_\reg/W)$ by mapping $T_i$ to the {\em inverse} 
of the Brieskorn loop $[\ga_i/W]\in \pi_1(\lt_\reg/W)$,
then we have established the theorem.
\epf

\begin{rmk}
In the later paper \cite{DP1565} a slightly different formula appears and here
we wish to clarify the (minor) discrepancy.
The action on $G^*$ descends along the map 
$\pi:G^*\to G^0; (b_-,b_+)\mapsto b_-^{-1}b_+$ to an action on the big cell
$G^0:=U_-TU_+\subset G$.
Corollary 14.4 on p.97 
of \cite{DP1565} gives the formula for this action on $G^0$ to be:
$$a=u_-t^2u_+\,\mapsto\, t^{-1}_i\xi_i(u_+)a\xi_i(u_+)^{-1}t_i.$$
Since $\pi$ is a 
covering map (corresponding to replacing $t=e^{\pi i\Lambda}$ by $t^2$) and
the action on $t$ is the standard Weyl group action, we deduce the
corresponding action on $G^*$ is as in \eqref{eqn: hol im in gstar3}, except
with each $t_i$ replaced by $t_i^{-1}$.
This action would be obtained from isomonodromy if we use the alternative 
construction of Tits' extended Weyl group noted at the end of Remark 
\ref{rmk: ewg}.
\end{rmk}

\end{section}

\begin{section}{Deformation of the Isomonodromy Hamiltonians} \label{sn: hams}
In this section the Hamiltonian description of the isomonodromic deformations
of the previous section will be given. 
From this the connection of De \!Concini--Millson--Toledano Laredo
will be
derived directly.

Let $\M^*:=\g^*\times\lt_\reg$ be the product of the dual of the Lie algebra
of $G$ with the regular subset of the chosen Cartan subalgebra.
View $\M^*$ as a trivial fibre bundle over $\lt_\reg$ with fibre $\g^*$.
Given $(B,A_0)\in\M^*$, consider the meromorphic connection
$A$ on the trivial $G$-bundle over $\IP^1$ associated to the 
$\g$-valued meromorphic one-form
\begin{equation}	\label{eqn: 1+2 conn}
A^s := \left(\frac{A_0}{z^2} + \frac{B}{z}\right)dz
\end{equation}
on $\IP^1$.
Restricting $A$ to the unit disc $\Delta$
(and using the compatible framing coming 
from the given trivialisation) specifies a point of the moduli
space $\M(A_0)$.
Thus there is a bundle map, 
$$\nu:\M^*\to\M;\qquad  (B,A_0)\mapsto A\vert_\Delta.$$
This map is holomorphic (by Corollary \ref{cor: hol dpce})
and it is easy to prove it is generically a local
analytic isomorphism. (It is studied fibrewise in Appendix \ref{last apx} and
in \cite{smpgafm}.)
Thus the isomonodromy problem for the connections  \eqref{eqn: 1+2 conn}
is essentially equivalent to that considered in the previous section.

The proofs of the following two lemmas are not significantly different 
from the $GL_n(\IC)$ case and so are
omitted here (cf. \cite{JMMS, Harn94, Dub95long, Hit95long}).

\begin{lem} \label{lem: imd eqns}
The pull-back along $\nu$ of the isomonodromy connection on $\M$
is given by the following non-linear differential equation
for sections $B:\lt_\reg\to\g^*$ of $\M^*$:
\begin{equation} \label{eqn: imds}
dB= \left[B, \ad^{-1}_{A_0}([dA_0,B])\right]
\end{equation}
where $d$ is the exterior derivative on $\lt_\reg$ and
$\g^*$ is identified with $\g$ via the Killing form.
(Note that $[dA_0,B]$ takes values in $\g^\od:= \bigoplus_{\al\in\R}\g_\al$
and that $\ad_{A_0}$ is invertible on $\g^\od$.)
\end{lem}

Clearly $B$ flows in a fixed coadjoint orbit in $\g^*$.
Thus, putting the standard Poisson structure on $\g^*$, one would expect
a symplectic interpretation. 
Indeed the equation \eqref{eqn: imds}  has the following time-dependent
Hamiltonian formulation.
Consider the one-form
$$\varpi:= \K\bigl(B,\ad^{-1}_{A_0}([dA_0,B])\bigr)$$
on $\M^*$, where $\K$ is the Killing form.
Given a vector field $v$ on $\lt_\reg$ there is a corresponding 
vector field $\wt v$ on $\M^*$ (zero in $\g^*$ directions)
and thus a function
$$H_v:= \langle \wt v,\varpi \rangle$$
on $\M^*$.

\begin{lem}
$H_v$ is a time-dependent Hamiltonian for the flow of the equation
\eqref{eqn: imds} along the vector field $v$
\end{lem}

\begin{rmk}
Usually one chooses a basis $\{ v_i\}$ of $\lt$ and writes
$\varpi = \sum H_i dt_i$ (where $A_0=\sum t_iv_i\in \lt$). 
It is this one-form which is used to define the isomonodromy
$\tau$ function.
\end{rmk}

Now observe that $\varpi$ may equivalently be viewed as a one-form 
on $\lt_\reg$ whose coefficients are quadratic polynomials on $\g^*$.
Let us identify these quadratic polynomials with $S^2\g=Sym^2\g$ and consider
the natural symmetrisation map $\phi:S\g\to U\g$ from the symmetric algebra to
the universal enveloping algebra.
(This corresponds to deforming the isomonodromy Hamiltonians under the
standard deformation `PBW quantisation' of $S\g$ into $U\g$.)

\begin{prop}
$$\phi(\varpi) = 
\sum_{\al\in \R_+} \frac{\K(\al,\al)}{2}(e_\al f_\al+f_\al e_\al)
\frac{d\al}{\al}\in U\g\otimes \Omega^1(\lt_\reg)$$
where $\{e_\al,f_\al, h_\al\}$ is the usual Chevalley basis for $\g$,
normalised so  that $[e_\al,f_\al]=h_\al$.

\end{prop}
The proof is a straightforward calculation.
This is precisely the $U\g$ valued one-form appearing in 
\cite{TL}; Given a representation $V$ of $\g$ and thus an algebra homomorphism 
$\rho:U\g\to \End(V)$, the flat connection of
De \!Concini--Millson--Toledano Laredo (whose holonomy is conjectured to give
the quantum Weyl group actions) is
$$d-h\rho(\phi(\varpi))$$ on the trivial vector bundle over $\lt_\reg$ with
fibre $V$, where $h\in \IC$ is constant.
\end{section}

\appendix
\begin{section}{Proofs for Section \ref{sn: sms}} \label{apx: pfs}

In the $GL_n(\IC)$ case these results appear in the paper 
\cite{BJL79} of Balser, Jurkat and Lutz, which in turn uses 
a theorem of Sibuya (that the map taking the Stokes multipliers is
surjective) and the main asymptotic existence theorem of Wasow
\cite{Was76} (in order to construct fundamental solutions).
Here we will follow the scheme of \cite{BJL79} wherever possible,
but notable exceptions arise in the use of both the above theorems:
1) The reduction to the asymptotic existence theorem is completely different
(see proof of Lemma \ref{lem: holom dpce}). Also an independent 
construction of fundamental solutions is given using `multisummation' rather
than the asymptotic existence theorem. 
2) For surjectivity, we instead follow the approach of Malgrange \cite{Mal79}
involving a $\overline \partial$-problem
which extends easily to general groups.

Note that we  must adapt the proofs from the $GL_n(\IC)$ case 
(rather than  simply choosing a faithful representation
$\g\hookrightarrow \gl_n(\IC)$ and using existing results)
because there are not 
representations taking elements of $\lt_\reg$ into regular diagonal
elements of $\gl_n(\IC)$ in general. 
\begin{eg} 
The standard representation of $\g=\so_4(\IC)$
is equivalent to writing
$$\g\cong\{ X\in\gl_4(\IC) \ \bigl\vert \ X^TJ + JX = 0 \}\subset \gl_4(\IC)$$
where 
$J=\left(
\begin{smallmatrix}
0 & I \\
I & 0 
\end{smallmatrix}
\right)$
and $I$ is the $2\times 2$ identity matrix.
This description is 
chosen so that we may take 
$\lt=\{ \diag(a,b,-a,-b) \ \bigl\vert\ a,b\in\IC \}$.
Now the regular elements of $\lt$ (as a Cartan subalgebra of $\g$) are
precisely those with both $a+b$ and $a-b$ nonzero.
However they will not be regular for $\gl_4(\IC)$ unless we also
impose $a\ne 0$ and $b\ne 0$ as well.
\end{eg}

One still may feel that for sufficiently generic values of the parameters 
one may always reduce to the $GL_n(\IC)$ case. 
Let us dispel this feeling:

\begin{lem} \label{lem: novlem}
There are reductive groups $G$ with the following property:
If $A_0\in\lt_\reg$ is any regular element of a Cartan subalgebra of
$\g=\Lia(G)$ and 
$\rho:\g\to\End(V)$ is any nontrivial representation of $\g$ then $\rho(A_0)$
does not have pairwise distinct eigenvalues.
\end{lem}
\pf
If $\rho(A_0)$ has pairwise distinct eigenvalues clearly each weight space of
$V$ is one-dimensional; $V$ is a multiplicity one representation.
However if for example $G=E_8$ then $\g$ has no nontrivial 
multiplicity one representations.
\epf

Thus there are groups for which one may {\em never} 
reduce to the $GL_n(\IC)$ case
via a representation.

For the purposes of this appendix we will use the notion of 
`Stokes directions' as well as the anti-Stokes directions already 
defined: $\si\in S^1$ is a {\em Stokes direction} if{f}
$\si-\pi/(2k-2)$ is an anti-Stokes direction.
(These are the directions along which the asymptotic behaviour of 
$\exp(\al\circ q)$ {\em changes} for some root $\al$, and they arise
as bounding directions of the supersectors.)
A crucial step enabling us to generalise \cite{BJL79} is the
following lemma. 

\begin{lem}	\label{lem: main apx lem}
Suppose $\theta\in S^1$ is not a Stokes direction.
Then the following hold

$\bullet$ The element
$$\lambda := \Real\bigl(A_0\exp(-i(k-1)\theta)\bigr)\in\lt_\IR$$
is in the interior of a Weyl chamber, and so determines an ordering
of the roots $\R$.

$\bullet$
The positive roots $\R^+(\lambda)$ (in this ordering) are precisely those 
roots supporting some anti-Stokes direction within $\pi/(2k-2)$ of $\theta$.

$\bullet$ Let $U_+$ be the unipotent part of the
Borel subgroup  $B_+\supset T$ determined by $\R^+(\lambda)$.
The following conditions on an element $C\in G$ are equivalent:

1) $z^\Lambda e^Q C e^{-Q}z^{-\Lambda}$ tends to $1\in G$ as $z\to 0$ in
the direction $\theta$.

2) $z^\Lambda e^Q C e^{-Q}z^{-\Lambda}$ is asymptotic to $1$ as $z\to
0$ in the direction $\theta$.   

3) $C\in U_+$.
\end{lem}

\pf
First we recall some group-theoretic facts (from e.g. \cite{Bor91}).
Let $X(T)=\Hom(T,\IC^*)$ be the character lattice of $T$, so that
$\R\subset X(T)$ naturally (by thinking of the roots multiplicatively).
In turn $\R$ is a subset of the real vector space 
$\lt_\IR^*:= X(T)\otimes_\IZ\IR$; This has (real) dual $\lt_\IR$ and
naturally $\lt_\IR\otimes_\IR\IC \cong \lt$.
By definition the Weyl chambers of $G$ relative to $T$ are the
connected components of 
\begin{equation}	\label{eq: weyl chambers}
\lt_{\IR,\reg}:= \{ \lambda\in\lt_\IR \ \bigl\vert\ 
\al(\lambda)\ne 0 \text{ for all } \al\in\R \}.
\end{equation}
Choosing a Weyl chamber is equivalent to choosing a system of positive
roots; If $\lambda\in\lt_{\IR,\reg}$ then the system of positive roots
corresponding to $\lambda$'s connected component is:
$$\R^+(\lambda) := \{ \al\in\R \ \bigl\vert \ \al(\lambda)>0 \}.$$

Now suppose
$\lambda := \Real\bigl(A_0\exp(-i(k-1)\theta)\bigr)$ as above and $\al\in\R$.
It is easy to check $\al(\lambda)=0$ if and only if
$\theta-\pi/(2k-2)$ is an anti-Stokes direction,
but by hypothesis this is not the case, so $\lambda$ is indeed regular.

Now consider the `sine-wave' function 
$f_\al(\phi):=\Real\bigl(\al(A_0)\exp(-i(k-1)\phi)\bigr)$ as
$\phi$ varies, for any $\al\in\R$. 
It has period $2\pi/(k-1)$ and is maximal at each 
anti-Stokes direction supporting $\al$.
Thus $f_\al(\theta)>0$ if{f} there is an anti-Stokes direction supported by
$\al$ within 
$\pi/(2k-2)$ of $\theta$.
In turn this is 
equivalent to $\al\in\R^+(\lambda)$, yielding the second statement.
(Note that, if $\arg(z)=\phi$, then $\Real(\al\circ q (z)) = -cf_\al(\phi)$
for some positive real $c$.)

For the third statement we will use the Bruhat decomposition of 
$G$ (cf. e.g. \cite{Bor91} 14.12).
Choose arbitrarily a lift $\wt w\in N(T)$ of each element 
$w\in W:= N(T)/T$ of the Weyl group.
The Bruhat decomposition says that $G$ is the disjoint 
union of the double cosets $B_+\wt w B_+$ as $w$ ranges over $W$.
The dense open `big cell' is the largest such coset (corresponding to the
`longest element' of $W$) and is equal to 
$B_-B_+$, where $B_-$ is the Borel subgroup opposite to $B_+$. 
Moreover the product map 
$$U_-\times T\times U_+\to G;\qquad (u_-,t,u_+)\mapsto u_-\cdot t\cdot u_+$$
is a diffeomorphism onto the big cell, where $U_\pm$ is the unipotent
part of $B_\pm$.

Now observe that each coset in the  Bruhat decomposition is stable under
conjugation by $T$, and that $U_+$ (and in particular the identity 
element of $G$) is in the big cell.
Thus we can reduce to the case where $C$ is in the big cell;
otherwise $3)$ is clearly not true and also neither of $1)$ or $2)$
hold,  since 
$z^\Lambda e^Q C e^{-Q}z^{-\Lambda}$ will remain outside of the big
cell.

Therefore if we label  
$\R^+(\lambda)=\{\al_1,\ldots \al_n\}$ and let $\al_{-i}=-\al_i$, 
then $C$ has a unique
decomposition 
$$C=u_{-1}\cdots u_{-n}tu_{1}\cdots u_n$$ 
with $u_{i}\in U_{\al_i}$ and $t\in T$ 
(\cite{Bor91} 14.5).
Each of these components is independent and 
$z^\Lambda e^Q C e^{-Q}z^{-\Lambda}= 
u^z_{-1}\cdots u^z_{-n}tu^z_{1}\cdots u^z_n$,
where $u^z_{i} = z^\Lambda e^Q u_{i} e^{-Q}z^{-\Lambda}\in
U_{\al_i}$.
Now, given a root $\al$ and $X\in\g_\al$, 
the key fact is that we know the behaviour of
$\Ad_{z^\Lambda e^Q} (X)$ as $z\to 0$ in the direction $\theta$; 
it decays exponentially if $\al\in\R^+(\lambda)$ and otherwise (if $X\ne 0$)
it explodes exponentially. 
(The dominant term of $z^\Lambda e^Q$ is $e^q$ and this acts on $X$ by
multiplication by $e^{\al\circ q}$, which has the said properties.)
Finally, in any representation $\rho$, $u_{i}$ is of the form
$1+\rho(X_i)$  with $X_i\in \g_{\al_i}$, so that
$C\in U_+$ if and only if   
$t=1$ and $X_{i}=0$ for all $i<0$, and
in turn (via the decomposition of  $z^\Lambda e^Q C
e^{-Q}z^{-\Lambda}$)
this is equivalent to both $1)$ and $2)$.
\epf

Now we will move onto the proofs of the results of Section
\ref{sn: sms}.

\pfms {\em(of Lemma \ref{lem: formal t}).}\ \ 
(This is adapted from 
\cite{JMU81} Proposition 2.2, 
\cite{Sab98} Theorem B.1.3 and
\cite{BV83} Lemma 1 p.42.)
For the existence of $\wh F$ and $A^0$ we proceed as follows. Write
$$A^s=A_0\frac{dz}{z^k}+  \cdots  +A_{k-1}\frac{dz}{z}+A_kdz+  \cdots$$
with $A_i\in\g$.  
First each $A_i$ will be moved into $\lt$ and then the nonsingular
part will be removed.
Let $\g^\od= \bigoplus_{\al\in\R}\g_\al$ (so that
$\ad_{A_0}:\g^\od \to \g^\od$ is an isomorphism), and 
let $\pr:\g\to\g^\od$ be the projection along $\lt$.
Suppose inductively that the first $p$ coefficients
$A_0,A_1,\ldots,A_{p-1}$ of $A^s$ are in $\lt$
(so the $p=1$ case holds by assumption).
By applying the gauge transformation $\exp(z^{p}H_{p})$ to $A^s$ 
(where 
$H_{p}\in \g$), we find
$$\exp(z^{p}H_{p})[A^s]= 
A^s+ [H_{p},A_0] z^{p-k}dz + O(z^{p-k+1})dz.$$ 
Thus $A_p+[H_{p},A_0]$ is the first coefficient which is not
necessarily in $\lt$,
and so by defining 
$$H_{p}:= (\ad_{A_0})^{-1}\bigl(\pr(A_p)\bigr)\in\g^\od$$
we ensure that the first $p+1$ coefficients of
$\exp(z^{p}H_{p})[A^s]$ are in $\lt$, completing the inductive step.
Hence if we define a formal transformation 
$\wh H\in\ G\flb z \frb$ to be the infinite product
$$\wh H :=\quad \cdots \exp(z^pH_p)\exp(z^{p-1}H_{p-1})\cdots\exp(zH_1)$$
then each coefficient of $\wh H[A^s]$ is in $\lt$.
Now define $A^0$ to be the principal part of $\wh H[A^s]$ so that
$\wh H[A^s]= A^0 + D$
with $D$ nonsingular.
Then define $\wt F :=e^{\left(-\int^z_0D\right)}\in T\flb z\frb$
(where
$\int^z_0 D \in \lt\flb z\frb$ is the series obtained from $D$ by replacing 
$z^p dz$ by $z^{p+1}/(p+1)$ for each $p\ge 0$),
so that 
$d\wt F(\wt F)^{-1}=d\log\wt F=-D$.
Thus $(\wt F\wh H)[A^s]=A^0$ and so 
$\wh F:=(\wt F\wh H)^{-1}\in G\flb z\frb$ 
is the desired formal transformation.

For the uniqueness it is clearly sufficient to show that if 
$\wh F[A^0]=A^1$, where $A^0$ and $A^1$ are formal types and $\wh F(0)=1$, 
then $\wh F=1$.
Now if $\wh F[A^0]=A^1$, 
it follows that $\wh F$ is actually convergent (since on using a
faithful representation we see $\wh F$ solves the diagonal system
$d\wh F = A^0 \wh F- \wh F A^1$). 
Let $F:\Delta\to G$ denote the sum of $\wh F$.
Also, if we write $A^i=dQ^i+\Lambda^idz/z$ for $i=0,1$, then
$d(e^{-Q^1}z^{-\Lambda^1} F z^{\Lambda^0}e^{Q^0}) = 0$, so that
$ F= z^{\Lambda^1}e^{Q^1} C e^{-Q^0}z^{-\Lambda^0}$ for some
constant $C\in G$.
Now $F\to 1$ on any sector and so (since $Q^0$ and $Q^1$ have the
same leading term) as in Lemma \ref{lem: main apx lem} we may deduce 
$C\in U_-\cap U_+=\{1\}$, and therefore $F=z^{\Lambda^1-\Lambda^0}e^{Q^1-Q^0}$.
The only way this can have a Taylor expansion with constant term 1 at $0$ is
if $\Lambda^1=\Lambda^0, Q^1=Q^0$ and so $\wh F = F =1$.
\epfms

\pfms {\em(of Lemma \ref{lem: s groups}).}\ \ 
Given a half-period $\bd\subset \IA$, let $\theta(\bd)$ be the
bisecting direction of the sector spanned by $\bd$.
By the symmetry of $\IA$, $\theta(\bd)-\pi/(2k-2)$ is
half-way between two consecutive anti-Stokes directions, 
so $\theta(\bd)$ is not a Stokes direction.
Therefore we may feed $\theta(\bd)$ into 
Lemma \ref{lem: main apx lem}, the second part of which immediately
yields the first statement of Lemma \ref{lem: s groups}.

The third statement of Lemma \ref{lem: s groups} 
is now immediate from Sections 14.5-14.8 of
\cite{Bor91}, (using the notion of `direct spanning' subgroups) and
then the second statement follows provided we check $\R(d)$ is a 
{\em closed} set of roots, in the sense that if 
$\al,\be\in\R(d)$ and $\al+\be\in\R$ then  $\al+\be\in\R(d)$.
This however is immediate from the definition of $\R(d)$.

For the fourth statement simply observe $\lambda$ is negated when
$\theta(\bd)$ is rotated by $\frac{\pi}{k-1}$.
\epfms

\pfms {\em(of Theorem \ref{thm: multisums}).}\ \ 

{\bf Uniqueness:} (cf. \cite{BJL79} Remark 1.4)
Suppose $F_1,F_2:\Sect_i\to G$ both have properties 
$1)$ and $2)$.
Thus $(F_1^{-1}F_2)[A^0]=A^0$ and so 
\begin{equation}	\label{eqn: Fratio}
F_1^{-1}F_2 = z^\Lambda e^Q C e^{-Q}z^{-\Lambda}
\end{equation}
for some constant $C\in G$. 
By $2)$,  $F_1^{-1}F_2$ extends to $\Ssect_i$ and is asymptotic to $1$
at zero there.
But $\Ssect_i$ has opening greater than $\pi/(k-1)$, so 
Lemma \ref{lem: main apx lem} implies $C\in U_+\cap U_-=\{ 1\}$ (by
taking two non-Stokes directions in $\Ssect_i$ differing by $\pi/(k-1)$).  
Hence $C=1$ and $F_1=F_2$.

{\bf Existence:}
Here we will use multisummation (see proof of Lemma \ref{lem: holom dpce}
for a more conventional approach).
Choose a faithful representation $G\hookrightarrow GL_n(\IC)$ such
that $T$ maps to the diagonal subgroup. 
Let $\ld\subset \gl_n(\IC)$ be the diagonal subalgebra and
let $\al_{ij}:\ld\to\IC; X\mapsto X_{ii}-X_{jj}$ be the roots of $GL_n(\IC)$.
Everything now will be written in this representation. 
Thus $\wh F$ is a formal solution to the system of linear differential
equations: 
\begin{equation}	\label{eqn: Feqn}
d\wh F = A\wh F-\wh FA^0.
\end{equation}
This equation has {\em levels}
$\bk := \{ -\deg(\al_{ij}\circ Q) \ \bigl\vert\  i,j=1,\ldots n 
\}\setminus\{0\}.$
Note that the highest level is $k-1$.
(If $k=2$ or if $G=GL_n(\IC)$ then this is 
the only level---however generally there may well be lower levels as well.)
The {\em singular directions} $\IA_\gl$
of \eqref{eqn: Feqn} are the $GL_n(\IC)$ 
anti-Stokes directions, defined as follows.
For each $i,j$, $\al_{ij}\circ Q$ is a polynomial in $1/z$ of
degree at most $k-1$. If $\al_{ij}\circ Q$ is not zero let 
$\IA^{ij}_\gl$ be the finite number of directions along which 
the leading term of $\al_{ij}\circ Q$ is real and negative
and let $\IA^{ij}_\gl$ be empty otherwise.
Then define $\IA_\gl$ to be the union of all these sets $\IA^{ij}_\gl$
as $i$ and $j$ vary.
(One may check that $\IA\subset \IA_\gl$.)
Then the main theorem in multisummation theory implies:
\begin{thm}[See \cite{BBRS91} Theorem 4.1] 	\label{thm: actual msums}
If $d\in S^1$ is not a singular direction then 
(each matrix entry of) $\wh F$ is 
$\bk$-summable in the direction $d$.
The $\bk$-sum of $\wh F$ along $d$ 
is holomorphic and asymptotic to $\wh F$ at zero in the sector
$\Sect(d-\frac{\pi}{2k-2}-\epsilon, d+\frac{\pi}{2k-2}+\epsilon)$
for some $\epsilon> 0$.
\end{thm}
Since they are unneeded here we omit discussion of the finer Gevrey
asymptotic properties that such sums possess,
although we do need the fact that multisummation is a morphism of differential
algebras (\cite{MR91} Theorem 1 p.348).
In more detail recall (\cite{BBRS91} Theorem 4.4) 
that the set $\IC\{z\}_{\bk,d}$
of formal power series which are $\bk$-summable in the direction $d$ 
is a differential subalgebra of $\IC\flb z \frb$.
Then $\bk$-summation maps this injectively onto some set
$\cO_{\bk,d}$ of germs at $0$ of holomorphic functions on 
$\Sect(d-\frac{\pi}{2k-2}, d+\frac{\pi}{2k-2})$.
Quite generally the map taking asymptotic expansions is easily seen to be a
differential algebra morphism, and so here it restrict to an 
{\em isomorphism} $\cO_{\bk,d} \mapright{\cong} \IC\{z\}_{\bk,d}$ 
of differential algebras.
By definition the multisummation operator is the inverse morphism.

Now to construct $\Sif$ choose any direction $d$ in $\Sect_i$ which is
not in (the finite set) $\IA_\gl$.
Let $\Sif$ be the multisum of
$\wh F$ along the direction $d$ from Theorem \ref{thm: actual msums}
and let $S$ be the sector appearing there.
Since multisummation is a morphism of differential algebras we deduce
first that $\Sif$ satisfies equation \eqref{eqn: Feqn} (as is standard in the
theory)
and secondly:

\begin{lem}
$\Sif$ takes values in $G$.
\end{lem}
\pf
This is because $G$, being reductive, 
is an affine algebraic group and so the matrix
entries of $\Sif$
satisfy the same polynomial equations as the entries of 
$\wh F$.
In more detail 
there are complex polynomials $\{p_j\}$ such that 
\begin{equation}	\label{eqn: G}
G \cong \{ (g,x)\in \IC^{n\times n}\times \IC \
\bigl\vert \ 
\det(g)\cdot x=1,\  p_j(g)=0 \ \forall \ j \},
\end{equation}
as a subgroup of $GL_n(\IC)$, for some $n$.
For any commutative algebra $R$ over $\IC$ 
the algebraic group $G(R)$ is defined
simply by replacing the two occurrences of $\IC$ in \eqref{eqn: G} by $R$.
Thus we wish to show $\Sif\in G(\cO(S))$ (the group of holomorphic
maps $S\to G$), given that 
$\wh F\in G\flb z\frb := G(\IC\flb z \frb)$.
But it is immediate that $p_j(\wh F)=0$ implies $p_j(\Sif)=0$ since 
multisummation is an algebra morphism. 
\epf

Next we must check that $\Sif$ has property $2)$ of Theorem 
\ref{thm: multisums}.  
The key point is that there are no Stokes directions in
$\Ssect_i\setminus S$; indeed the Stokes directions in the supersector
$\Ssect_i$ closest to the boundary rays are 
$d_{i+1}-\frac{\pi}{2k-2}$ and $d_{i}+\frac{\pi}{2k-2}$, both of 
which are in $S$.
Thus the following $G$-valued analogue of the extension lemma of
\cite{BJL79} will yield $2)$:
\begin{lem}[cf. \cite{BJL79} Lemma 1 p.73] \label{lem: extn}
Suppose $S,\wt S$ are two sectors with non-empty intersection and such
that $\wt S$ contains no Stokes directions.
If $F: S\to G$ is a holomorphic map asymptotic to $\wh F$ at $0$ in
$S$ and such that $F[A^0]=A$, then the analytic continuation of $F$
to $S\cup\wt S$ is asymptotic to $\wh F$ at $0$ in $S\cup\wt S$.
\end{lem}
\pf
Choose any holomorphic map $\wt F:\wt S\to G$ asymptotic to $\wh F$ at
$0$ and such that $F[A^0]=A$ (using multisummation for example---the
hypotheses imply $\wt S$ has opening $<\pi/(k-1)$).
Then (as in the uniqueness part above) 
there exists a constant $C\in G$ such that 
$F=\wt F z^\Lambda e^Q C e^{-Q}z^{-\Lambda}$ in $S\cap\wt S$.
Thus $\wt F z^\Lambda e^Q C e^{-Q}z^{-\Lambda}$ is the analytic
continuation of $F$ to $\wt S$.
Now since $\wt F^{-1}F$ is asymptotic to $1$ at $0$ in $S\cap\wt S$,
Lemma \ref{lem: main apx lem} implies $C\in U_+$, where the root
ordering is determined by any $\theta$ in  $S\cap\wt S$.
But since $\wt S$ contains no Stokes directions, 
Lemma \ref{lem: main apx lem} implies 
$z^\Lambda e^Q C e^{-Q}z^{-\Lambda}$ 
is asymptotic to $1$ at $0$ in all of $\wt S$.
In turn it follows that the analytic continuation of $F$ is asymptotic
to $\wh F$ on all $S\cup \wt S$.
\epf

Finally the last
statement of Theorem 
\ref{thm: multisums} is immediate either from the morphism properties of
multisummation, or from uniqueness.
\epfms

\pfms {\em(of Lemma \ref{lem: sf in sfg}).}\ \ 
To see $S_j\in\ISto_\bd(A^0)$ recall from Lemma \ref{lem: s groups}
that $\ISto_\bd(A^0)=U_+$ where the
root order is determined by the bisecting direction $\theta(\bd)$ of 
$\bd$.
Now observe $\Ssect_{(j-1)l}\cap\Ssect_{jl}$ contains $\theta(\bd)$
and so by Theorem \ref{thm: multisums} (if $j\ne 1$)
$z^\Lambda e^Q S_j e^{-Q}z^{-\Lambda} = 
\Sigma_{jl}(\wh F)^{-1}\Sigma_{(j-1)l}(\wh F)$ is asymptotic to $1$
along $\theta(\bd)$.
Thus Lemma \ref{lem: main apx lem} implies $S_j\in\ISto_\bd(A^0)$.
(For $j=1$ the argument is the same once the change in branch of
$\log(z)$ is accounted for.)
In turn
to see $K_i\in \ISto_{d_i}(A^0)$ simply observe
$\ISto_{d_i}(A^0) = \ISto_\bd(A^0)\cap\ISto_\bdp(A^0)$
(where $\bd=(d_i,\ldots,d_{i+l-1})$ and 
$\bdp=(d_{i-l+1},\ldots,d_i)$ are the two half-periods ending on $d_i$),
and that the above argument implies $K_i$ is in this intersection,
since $\Sigma_{i}(\wh F)^{-1}\Sigma_{i-1}(\wh F)$ is asymptotic to $1$
along both $\theta(\bd)$ and $\theta(\bdp)$.
\epfms

Now we will 
establish the surjectivity of the irregular Riemann--Hilbert map in
Theorem \ref{thm: irr RH im}.
Fix a formal type $A^0$ and let $\IA$ be the corresponding set of anti-Stokes
directions.
Also fix a choice of initial sector and branch of $\log(z)$ as in Section
\ref{sn: sms}.
Now choose arbitrarily a Stokes factor $K_d\in\ISto_d(A^0)$ for each
$d\in \IA$.

\begin{thm}
There exists a meromorphic connection $A$ on the trivial principal $G$-bundle
over $\Delta$ having formal type $A^0$ and Stokes factors 
$\{K_d\}$.
\end{thm}

\pf
(This is an adaptation of \cite{BalserBk} Section 9.7, 
except we replace the key step with a $\overline\partial$-problem,
as was suggested by Malgrange \cite{Mal79} 
and fleshed out in \cite{BV89} Section 4.4.)
First we remark that it is sufficient to construct $A$ only in a 
neighbourhood of the origin 
because any such connection is gauge equivalent to a connection
defined over the whole disc.
(One may prove this as follows:
Given $A$ over $\Delta_\epsilon:=\{z\ \bigl\vert\ \vert z\vert\le\epsilon\}$,
choose any holomorphic connection $A^1$ on $G\times \Delta^*$ with the same
monodromy as $A$ around $0$.
The ratio $\Phi^1\cdot\Phi^{-1}$ of corresponding fundamental solutions then 
defines a holomorphic map from $\Delta^*_\epsilon$ to $G$
which we 
use as a clutching function to define a principal $G$-bundle $P$ over 
$\Delta$. The connections $A,A^1$ define a single meromorphic
connection on $P$.
Moreover $P$ is trivial since all $G$-bundles over a disk are---cf. 
\cite{GrauertSP} p.370.)

Now choose $j$ such that $K_{d_j}\ne 1$.
By induction on the number of non-trivial Stokes factors we may assume there
is a connection $B$ having Stokes factor $K_d$ for each $d\ne d_j$ but having
Stokes factor $1\in \ISto_{d_j}(A^0)$ along $d_j$.
Write $K= K_{d_j}$ for simplicity and let $B^s=-s^*(B)$ as usual.
(Here $s$ is the identity section of the trivial $G$-bundle $G\times\Delta$.)
Let $\chi_j(z)$ be the canonical fundamental solution of $B$ on $\Sect_j$ from
Definition \ref{def: sf2} and define 
$\chi(z):= \chi_j K \chi_j^{-1}$.
Let $\wt d_j$ be the lift of the direction $d_j$ to the universal cover
$\wt \Delta^*$ of the punctured disk determined by the chosen branch of
$\log(z)$ and let $\wt\Sect_j$ be the lift of $\Sect_j$.
Let $\cS^*$ denote the sector in 
$\wt \Delta^*$ (of opening more than $2\pi$) 
from $\wt d_j-\frac{\pi}{2k-2}+\delta$ to
$\wt d_j + 2\pi +\frac{\pi}{2k-2}-\delta$,
and let $\cS = \Sect(\wt d_j-\frac{\pi}{2k-2}+\delta, 
\wt d_j +\frac{\pi}{2k-2}-\delta)$.
Here $\delta>0$ is fixed so that no Stokes directions lie in the interior
of either component of $\Sect(\wt d_j-\frac{\pi}{2k-2}, 
\wt d_j +\frac{\pi}{2k-2})\setminus \cS$ 
(i.e. $\delta < \min\{\vert d_j-d_{j\pm 1}\vert\}$).
We now claim that
there exists a holomorphic map
$\tau:\cS^*\to G$ having an asymptotic expansion with constant term
$1$ in 
$\cS^*$ and such that
\begin{equation}	\label{eqn: tau prop}
\tau(ze^{2\pi i})=\tau(z)\chi(z)
\end{equation}
for any $z\in \wt \Sect_j\cap \cS$, where $\chi$ is pulled up to
$\wt \Sect_j$ in the obvious way.

Such $\tau$ may be constructed as follows.
Let $\cS' = \Sect(\wt d_j-\frac{\pi}{2k-2}+\delta, \wt d_j+\pi)$
and (as in \cite{BV89}  Lemma 4.3.2 and using the
exponential map for $G$)
extend $\chi$ to a $C^\infty$ map $f:\cS'\to G$ such that 
$f\vert_\cS = \chi$,  
$f(z)= 1$ for $\arg(z)$ in some neighbourhood
of $\wt d_j+\pi$ and $f\sim 1$ on all $\cS'$.
(By construction $\chi\sim 1$ on $\Sect(\wt d_j-\frac{\pi}{2k-2}, 
\wt d_j +\frac{\pi}{2k-2})$.)
Then define a $C^\infty$ $\g$-valued one-form $\al$ on $\Delta$ by
letting $\al= f^{-1}\overline\partial f $ on $\cS'$ and extending by zero.
Now solve the $\overline\partial$-problem
$g^{-1} \overline\partial g  = \al$ for a smooth map $g$ from some
neighbourhood of $0\in \Delta$ to $G$, with $g(0)=1$.
(This is possible for the same reasons as in the $GL_n(\IC)$ case,
for which cf. e.g. \cite{AB82} p.555.)
Finally define $\tau: \cS^*\to G$ by $\tau= gf^{-1}$ for
$\arg(z) \le \wt d_j +\pi$ and $\tau= g $ for
$\arg(z) \ge \wt d_j +\pi$; one easily checks this has the properties claimed.

To complete the proof 
define $\wt \chi(z):= \tau(z)\chi_j(z)$ for $z$ in 
$\cS^*$ (where $\chi_j$ is continued from $\wt\Sect_j$ as a 
fundamental solution of $B$).
Then \eqref{eqn: tau prop} implies 
$A^s:=(d\wt\chi)\wt\chi^{-1}$ is invariant under rotation by $2 \pi$ and
so defines a 
$\g$-valued one-form on a neighbourhood of $0$ in $\Delta^*$.
We will show that the corresponding connection $A$ on the trivial principal
$G$-bundle has the desired properties.
First observe that $A^s= \tau[B^s]$ by holomorphicity, since this certainly
holds near $0$ in $\wt \Sect_j$.
Since $\tau$ admits an asymptotic expansion 
$\wh \tau\in G\flb z\frb$ in $\cS^*$
it follows that $A^s$ admits {\em Laurent} expansion $\wh\tau[B^s]$.
Thus if $B^s=\wh F[A^0]$ (from Lemma \ref{lem: formal t}) then
$A^s = (\wh \tau\circ\wh F)[A^0]$, and so $A$ has formal type $A^0$ as 
required.
Now from the range of validity of the asymptotic expansion of $\tau$, and from
the uniqueness of the sums in Theorem \ref{thm: multisums},
we deduce that for $z\in\Sect_i$:
$$\Sigma_i(\wh \tau\circ\wh F)(z) = \tau(\wt z)\cdot\Sigma_i(\wh F)(z)$$
where $\wt z\in \wt \Delta^*$ lies over $z$ and between directions
$\wt d_j$ and $\wt d_j+2\pi$.
(On $\Sect_{j-1}$ and $\Sect_j$ one needs to use the extension lemma,
Lemma \ref{lem: extn}, 
as well---which is applicable by the choice of $\delta$.)
In turn we immediately find that $A$ has the same Stokes factors as $B$ except
in $\ISto_{d_j}(A^0)$. 
Here (across $d_j$) by definition $A$ has Stokes factor
$$\Phi_j^{-1}\Phi_{j-1} = 
\chi_j^{-1}\tau^{-1}(z)\tau(ze^{2 \pi i})\chi_{j-1}$$
if $j\ne\#\IA$,
where $\Phi_i, \chi_i$ denote canonical solutions of $A,B$ respectively.
Since $B$ has trivial Stokes factor here   $\chi_{j-1}=\chi_j$ and so
by \eqref{eqn: tau prop} and the definition of $\chi$ we find
$\Phi_j^{-1}\Phi_{j-1} = \chi_j^{-1}\chi\chi_{j} = K$
as required. (Similarly if $j=\#\IA$.)
\epf

\pfms {\em(of Lemma \ref{lem: holom dpce}).}\ \ 
Here an alternative construction of the `sums' of Theorem \ref{thm: multisums}
will be given, closer to the usual approach in the $GL_n(\IC)$ case.
This is more direct
than the multisummation approach above and enables us to
prove that the sums vary holomorphically with parameters.

The usual construction of $\Si_0(\wh F)$ (cf. e.g. \cite{Was76}) consists of
two steps.
Roughly speaking 
one converts the equation satisfied by $\wh F$ into two independent
nonlinear equations for `parts' of $\wh F$.
Then an asymptotic existence theorem is used to find analytic solutions to
these two equations, that are asymptotic to the corresponding 
parts of $\wh F$.   
Usually for $G=GL_n(\IC)$ (see \cite{Was76} \S 12.1, \cite{Sib62}) 
the two
equations involve the upper and lower triangular parts of $\wh F-1$.
For 
general $G$ this makes no sense: an alternative procedure will be used here
to reduce the problem to the {\em same} asymptotic existence theorem.

If $a$ and $b$ are the boundary directions of the sector $\check S$ (so
$\check S=\Sect(a,b)$), let 
$S:=\Sect\left(b-\frac{\pi}{2k-2}, a+\frac{\pi}{2k-2}\right)$
(a sector of opening less than $\frac{\pi}{k-1}$ centred on $\check S$).
For simplicity write $F=\Si_0(\wh F)$ for the $G$-valued map we are seeking. 
This should have asymptotic expansion $\wh F$ as $z\to 0$ in $S$ and should
solve (for each $x\in U$) the equation
\begin{equation} \label{eqn: Feqn2}
(dF)F^{-1} = A-FA^0F^{-1} 
\end{equation}
on the sector $S\subset \Delta$, where $d$ is the exterior derivative on
$\Delta$. (Here $(dF)F^{-1}$ is defined in the usual way as the pullback of
the right-invariant Maurer--Cartan form on $G$ under the map $F$
and $FA^0F^{-1}:=\Ad_FA^0$.)
Note that such $F$ is unique since (for each $x$) the extension lemma
(Lemma \ref{lem: extn}) says $F$ extends to $\Ssect_0$ maintaining the
asymptotic expansion $\wh F$, and so the uniqueness part of Theorem
\ref{thm: multisums} fixes $F$.

Now, because the big cell $G^0=U_-TU_+\subset G$ is open and contains the
identity, we find
\begin{lem} \label{lem: factns}
(1) $\wh F$ admits a unique factorisation 
$$\wh F= \wh u_-\cdot \wh t\cdot \wh u_+$$
with $\wh u_\pm\in U_\pm\bigl(\cO(X)\flb z\frb\bigr)$ and 
$\wh t\in T\bigl(\cO(X)\flb z\frb\bigr)$.

(2) For $z\in S$ sufficiently close to $0$, any solution $F$ of 
\eqref{eqn: Feqn2}  asymptotic to $\wh F$ has a unique factorisation
$$ F= u_- \cdot t \cdot u_+$$
with $u_\pm,t$ taking values in $U_\pm, T$ respectively.
\end{lem}

Substituting $F= u_- \cdot t \cdot u_+$ into \eqref{eqn: Feqn2} and rearranging
yields:
\begin{equation} \label{eqn: fulleqn}
u_-^{-1}du_- + (dt)t^{-1} + t(du_+)u_+^{-1}t^{-1} + tu_+A^0u^{-1}_+t^{-1} -
u_-^{-1}Au_- =0.
\end{equation}
Taking the $\lu_-$ component of this gives the independent equation
\begin{equation} \label{eqn: -eqn}
u_-^{-1}du_- = \pi_-(u_-^{-1}Au_-)
\end{equation}
for $u_-$, where $\pi_-:\g=\lu_-\oplus\lt\oplus\lu_+\to\lu_-$ is the
projection.  We wish to solve \eqref{eqn: -eqn} using the following asymptotic
existence theorem with parameters, which is also used in the $GL_n(\IC)$ case.

\begin{thm} \label{thm: aep}
Let $S$ be an open sector in the complex $z$ plane with vertex $0$ and opening
not exceeding $\frac{\pi}{k-1}$.
Let $f(z,L,x):\Delta\times\IC^N\times U\to \IC^N$ be a holomorphic map such
that

(i) The Jacobian matrix 
$\left(\frac{\partial f_i}{\partial(L)_j}\right)\Bigl\vert_{L=0,z=0}$ 
is invertible for all $x\in \overline U$, and

(ii) The differential equation
\begin{equation} \label{eqn: nlin feqn}
\frac{dL}{dz} = \frac{f(z,L,x)}{z^k}
\end{equation}
admits a formal power series solution $\wh L= \sum_1^\infty L_r(x)z^r\in 
\IC^N\flb z\frb\otimes \cO(U)$.

Then there exists, for sufficiently small $z\in S$, a holomorphic solution
$L(z,x)$ of \eqref{eqn: nlin feqn} having asymptotic expansion 
$\wh L$ in $S$ uniformly in some neighbourhood of $x_0\in U$.
\end{thm}
\pf
Without parameters this is Theorem 14.1 of \cite{Was76}.
The method of successive approximations used there extends immediately to
the case with parameters since the uniform convergence of successive
approximations does not destroy holomorphicity with respect to parameters
(cf. \cite{Sib71} Remark 2 p.161). 
Alternatively Sibuya proves a very similar result (\cite{Sib62} Lemma 2)
for the case where the parameters become singular on a sector. 
As remarked in \cite{HsSib66} p.100 the case when the parameters are
nonsingular and on a disc (as required here)
is proved in exactly the same manner. 
\epf

The trick to convert \eqref{eqn: -eqn} into the form \eqref{eqn: nlin feqn}
is as follows. (We will have $\IC^N=\lu_-$.)
Since $U_-$ is unipotent the exponential map $\exp:\lu_-\to U_-$ is an
algebraic isomorphism, and its derivative gives an isomorphism
$\exp_*:T\lu_-\to TU_-$ of the tangent bundles.
If we identify $T\lu_-\cong \lu_-\times\lu_-$ using the vector space structure
of $\lu_-$ and $TU_-\cong U_-\times\lu_-$ using left multiplication in $U_-$
then we deduce:
\begin{lem} \label{lem: psi im}
There is an algebraic isomorphism 
$$\psi:\lu_-\times\lu_- \mapright{\cong} U_-\times\lu_-$$
which is linear in the second component and such that
if $L(z):\Delta\to \lu_-$ is any holomorphic map then
$$\psi\left(L,\frac{dL}{dz}\right) = \left(e^L, e^{-L}\frac{d}{dz}(e^L)\right).$$
\end{lem}

Thus, setting $u_-=e^L$, equation \eqref{eqn: -eqn} is equivalent to
\eqref{eqn: nlin feqn} with the map $f$ defined by
$$(L,f(z,L,x)) = \psi^{-1}(u_-, h(z,u_-,x))\in\lu_-\times\lu_-$$
for any $z\in \Delta, L\in \lu_-,x\in U$, where
$u_-=e^L$ and 
$h(z,u_-,x):= \langle z^k\pi_-(u_-^{-1}A(z,x)u_-), 
\frac{\partial}{\partial z}\rangle$ is from 
\eqref{eqn: -eqn}.

Clearly the formal solution $\wh u_-$ of \eqref{eqn: -eqn} induces a 
formal solution of \eqref{eqn: nlin feqn} of the desired form and so all that
remains to solve \eqref{eqn: nlin feqn} is to check the Jacobian condition
in Theorem \ref{thm: aep}.
Geometrically this condition says precisely that the graph
$\Ga(f)\subset \lu^h_-\times\lu^v_-$ of the map
$f(0,\cdot,x):\lu^h_-\to\lu^v_-$ is transverse to the horizontal subspace
$\lu_-^h$ at $L=0$.
(Here $\lu^h_-,\lu^v_-$ are just copies of $\lu_-$ labeled `horizontal' and
`vertical'; $f$ is viewed as a section of the tangent bundle to
$\lu_-^h$.)
Since $\psi$ is a diffeomorphism it is sufficient to check this transversality
condition on $U_-\times\lu_-$.
By definition the graph of $f$ corresponds to the graph of $h$ under $\psi$.
Clearly $L=0$ if{f} $u_-=1$   and $h(1) = \pi_-(A_0)=0$ since $A_0\in\lt$.
(Here we omit the arguments $z=0$ and $x$ of $h$ for notational simplicity.) 
Thus $\psi(0,f(0)) = (1,0)$ and we must check that the two vector spaces
$\psi_*(\lu_-\times\{0\})$ and $T_{(1,0)}\Ga(h)$ are transverse subspaces 
of the tangent space $T_{(1,0)}(U_-\times\lu_-) = \lu_-\times\lu_-$.
The tangent space to the graph of $h$ at $(1,0)$ is
$$\{(X,\pi_-[X,A_0]) \ \bigl\vert \ X\in\lu_- \}$$
and one may calculate that the derivative of $\psi$ maps the horizontal
subspace $\lu_-\times\{0\}$ to 
$$\{(X,0) \ \bigl\vert \ X\in\lu_- \}.$$
Since $A_0(x)\in\lt_\reg$ it follows immediately 
that these two subspaces are indeed transverse.
Thus we may apply Theorem \ref{thm: aep} to obtain a holomorphic 
solution $L(z,x)$ of \eqref{eqn: nlin feqn} and in turn obtain a solution 
$u_-=e^L$ of \eqref{eqn: -eqn}.

Given this solution $u_-$, now consider the $\lt$ component 
\begin{equation} \label{eqn: t cmpt}
(dt)t^{-1} = \delta(u_-^{-1}Au_-) - A^0
\end{equation}
of the full equation \eqref{eqn: fulleqn}, where 
$\delta:\g=\lu_-\oplus\lt\oplus\lu_+\to \lt$ is the projection.
This equation has formal solution $\wh t$ and so the right-hand
side of  \eqref{eqn: t cmpt} has {\em nonsingular} asymptotic expansion as
$z\to 0$ in $S$.
Immediately this implies \eqref{eqn: t cmpt} has a unique holomorphic solution
tending to $1\in T$ as $z\to 0$, given by
$$t(z,x) = \exp\left(\int_0^z \bigl(\delta(u_-^{-1}Au_-) - A^0\bigr)\right)$$
(cf. \cite{Was76} Theorem 8.7 p.38).

Thus we have obtained all except the $u_+$ component 
of the desired solution $F=u_-tu_+$.
To obtain $u_+$ we repeat all the above procedure with the opposite
factorisation $F=w_+sw_-$ of $F$ (with $w_\pm\in U_\pm$ and $s\in T$).
This yields $w_+$ and $s$.
Then, for sufficiently small $z$, the components $u_+$ and $w_-$ are determined
(holomorphically) from $u_-,w_+,s,t$ by the equation
$$w_+^{-1}u_-t=sw_-u_+^{-1}$$
since both Bruhat decompositions are unique and the left-hand side is known.
The resulting solution $F=u_-tu_+=w_+sw_-$ then has the desired properties.
\epfms

\end{section}

\begin{section}{} \label{last apx}

In this appendix we will explain how the results of Section 
\ref{sn: sms} enable us to extend Theorems 1 and 2 of 
\cite{smpgafm} from $GL_n(\IC)$ to arbitrary connected complex 
reductive groups $G$.
The main modifications of the proofs in \cite{smpgafm} are purely notational
and so here we will concentrate on giving a clear statement of the results.
The set-up is as follows.

Let $K$ be any compact connected Lie group.
Choose a maximal torus $T_K\subset K$ and
a non-degenerate symmetric invariant bilinear form $\K$ on
$\lk=\Lia(K)$.
(Thus if $K$ is semisimple we may take $\K$ to be the Killing form, or if
$K=U(n)$ then $\K(A,B)=\tr(AB)$ will do.)

Let $G$ be the complex algebraic group associated to $K$ 
(as in \cite{Chevbook}; $G$ is the variety associated to the 
complex representative ring of $K$).
Any complex connected reductive Lie group $G$ arises in this way 
(see \cite{HochsStr}).
We have $\g=\lk\otimes \IC$ with Cartan subalgebra
$\lt = \lt_K\otimes\IC$ and $G$ has maximal torus
$T=\exp(\lt)$.
Extend $\K$ $\IC$-bilinearly to $\K:\g\otimes\g\to\IC$. 
The group $G$ comes equipped with an involution (with fixed point set
canonically isomorphic to $K$), which we will denote by
$g\mapsto g^{-\dagger}$. (It is denoted $\iota$ in \cite{Chevbook}.)
This induces an anti-holomorphic involution of $\g$ 
(to be denoted $X\mapsto -X^\dagger$) fixing $\lk$ pointwise.

Note that $\lt$ comes with two real structures: one ($A\mapsto -A^\dagger$)
from the
identification $\lt = \lt_K\otimes\IC$ and another 
(to be denoted $A\mapsto \overline A$) defined via
the identification $\lt = \lt_\IR\otimes\IC$ where
$\lt_\IR:= X_*(T)\otimes_\IZ \IR$. 
(Here $X_*(T):= \Hom(\IC^*,T)$ is embedded in $\lt$ by differentiation.)
One may check $\lt_K=i\lt_\IR$.

Apart from the choices $K, \lt_K, \K$ made so far we need to make three
further
choices in order to define the monodromy map 
\begin{equation} \label{eqn: monmap}
\nu:\g^*\longrightarrow G^*.
\end{equation}
These are:

1) A regular element $A_0\in\lt_\reg$. This determines anti-Stokes directions
   etc. as in Section \ref{sn: sms} (taking the pole order $k=2$).

2) An initial sector $\Sect_0$ bounded by two consecutive anti-Stokes
   directions, and

3) A choice of branch of $\log(z)$ on $\Sect_0$.

The choice of initial sector determines a system of positive roots
$\R(d_1)\cup\cdots\cup\R(d_l)$
as in Lemma \ref{lem: s groups} (with $l=\#\IA/2$) and thus a Borel subgroup
$B_+\subset G$ containing $T$.
Let $B_-$ be the opposite Borel subgroup and define the dual Poisson Lie group 
$G^*$ as in \eqref{gstar} in terms of $B_\pm$.
$G^*$ is a contractible Lie group of the same dimension as $G$ and has
a natural Poisson Lie group structure which may be defined directly and 
geometrically as for $GL_n(\IC)$ in \cite{smpgafm}.

The monodromy map \eqref{eqn: monmap} is then defined as follows.
Given $B\in \g^*$, consider the meromorphic connection on the trivial
principal $G$-bundle over the unit disc $\Delta$ determined by the 
$\g$-valued meromorphic one-form
$$A^s := \left(\frac{A_0}{z^2} + \frac{B}{z}\right)dz,$$
where $B$ is viewed as an element of $\g$ via $\K$.
This connection has Stokes multipliers $(S_+,S_-)=(S_1,S_2)\in U_+\times U_-$
defined in Definition \ref{def: sf2} above, using choices 2) and 3),
and so determines an element
$$(b_-,b_+,\Lambda)\in G^*$$
via the formulae:
$$b_- = e^{-\pi i \Lambda} S_-^{-1}, \quad 
  b_+ = e^{-\pi i \Lambda} S_+e^{2\pi i\Lambda},\quad
  \Lambda= \delta(B)$$ 
so that $b_-^{-1}b_+ = S_-S_+\exp(2\pi i\Lambda)$.
(Here $\delta:\g\to \lt$ is the projection with kernel
$\lu_+\oplus\lu_-$.)
A slightly more direct/elegant definition of $\nu$ may be given, 
without first going through Stokes multipliers, exactly 
as before in Section 4 of \cite{smpgafm}.

The monodromy map $\nu$ is a holomorphic map by Corollary \ref{cor: hol dpce}
and it is easy to prove (as in \cite{smpgafm}) it is generically a 
local analytic isomorphism and any generic symplectic leaf of 
$\g^*$ maps into a symplectic leaf of $G^*$.
The approach of \cite{smpgafm} extends to yield:
\begin{thm} \label{thm: poissonness}
The monodromy map $\nu$ is a 
Poisson map for each choice of $A_0, \Sect_0, \log(z)$,
where $\g^*$ has its standard complex Poisson structure and 
$G^*$ has its canonical complex Poisson Lie group structure, but
scaled by a factor of $2\pi i$.
\end{thm}
\pf
As in  \cite{smpgafm};
Just replace any expression of the form $\tr(AB)$ by $\K(A,B)$
and any reference to the {\em difference of eigenvalues} of any element
$A\in \g$ (which now makes no sense), by 
the {\em eigenvalues} of $\ad_A\in \End(\g)$. 
(The left and right-invariant Maurer--Cartan forms on $G$ make sense of
expressions of the form $H^{-1}H'$ and $H'H^{-1}$ for maps 
$H$ into $G$.) 
The only minor subtlety 
is in the proof of \cite{smpgafm} Lemma 27 where one needs
the fact that $\Ad_{(e^Q)}(n_-)$ tends to zero as $z\to 0$ along a certain 
direction
($-\theta$) for any fixed $n_-\in \lu_-$. However this follows 
directly from the third part of Lemma \ref{lem: main apx lem} of the present
paper.
\epf

\begin{rmk}
Thus locally the monodromy maps give appropriate `canonical' 
coordinate changes to
integrate the explicit non-linear isomonodromy equations \eqref{eqn: imds}.
This indicates just how complicated the monodromy maps are: 
for $G=SL_3(\IC)$ equation \eqref{eqn: imds} is equivalent to the full family
of Painlev\'e VI equations---generic solutions of which are known to involve 
`new'
transcendental functions.
\end{rmk}

Now suppose $A_0$ is purely imaginary ($\overline A_0=-A_0$).
Then there are only two anti-Stokes directions; the two halves of the
imaginary axis. Take $\Sect_0$ to be the sector containing the
positive real axis $\IR_+$ and use the branch of $\log(z)$ 
which is real on $\IR_+$. 
One may then check (as in \cite{smpgafm} Lemma 29) that if
$(b_-,b_+,\Lambda)=\nu(B)$ then
$\nu(-B^\dagger)=(b_+^{-\dagger},b_-^{-\dagger},-\Lambda^\dagger)$
so that $\nu$ restricts to a (real analytic) map
$$\nu\vert_{\lk^*}:\lk^*\mapright{} K^*,$$
where $\lk^*\cong\lk$ via $\K$ and $K^*\subset G^*$
is defined to be the fixed point subgroup of the involution
$$(b_-,b_+,\Lambda)\mapsto(b_+^{-\dagger},b_-^{-\dagger},-\Lambda^\dagger).$$
The group $K^*$ has a natural (real) Poisson Lie group structure which may be
defined as in \cite{smpgafm} for $K=U(n)$.
All of these restricted monodromy maps are Ginzburg--Weinstein isomorphisms:
\begin{thm}	\label{thm: GW ims}
If $A_0$ is purely imaginary then the corresponding 
monodromy map restricts to a (real) Poisson diffeomorphism
$\lk^*\cong K^*$ from the dual of the Lie algebra of $K$ to the
dual Poisson Lie group 
(with its standard Poisson structure, scaled by a factor of $\pi$).
\end{thm}
\pf
The proof in \cite{smpgafm} goes through once the notational changes 
given in the previous proof are made.
The fact that the `unique Hermitian logarithms' appearing in the proof of
Lemma 31 \cite{smpgafm} still exist 
(Hermitian now meaning $-X^\dagger=-X\in \g$),
follows easily from the fact that $G$ has a faithful representation 
$\rho:G\hookrightarrow GL_N(\IC)$ with $\rho(K)=\rho(G)\cap U(N)$
(cf. \cite{Chevbook} Lemma 2 p.201). 
\epf

\begin{rmk}
1) The permutation matrices used in \cite{smpgafm} have now been banished;
Consequently the group $G^*$ now depends on the choice of initial sector
(a priori we make no choice of positive roots). The pleasant effect is
that $\nu$ is now always $T$-equivariant, with $T$ acting on $\g^*$ via the
coadjoint action and on $G^*$ via the left or right dressing action. 
(The left and right dressing actions agree when restricted to $T$.)
In turn the Ginzburg--Weinstein isomorphisms constructed above are all
$T_K$-equivariant.

2) The new proof of the theorem of Duistermaat given in Section 6 of
\cite{smpgafm} for $GL_n(\IC)$
also extends immediately to connected complex reductive $G$. 
\end{rmk}

\end{section}

\renewcommand{\baselinestretch}{1}		%
\normalsize
\bibliographystyle{amsplain}	\label{biby}
\bibliography{../thesis/syr}	
\end{document}